\documentclass[11pt]{article}
\usepackage{amsmath}
\usepackage{amssymb}
\usepackage{amscd}
\usepackage{verbatim}
\newcommand{\be}{\begin{equation}}
\newcommand{\bs}{\begin{sub}}
\newcommand{\es}{\end{sub}}
\newcommand{\bsn}{\begin{subn}}
\newcommand{\esn}{\end{subn}}
\newcommand{\bea}{\begin{eqnarray}}
\newcommand{\eea}{\end{eqnarray}}
\newcommand{\BA}[1]{\begin{array}{#1}}
\newcommand{\EA}{\end{array}}
\newcommand{\Real}{\mathbb{R}}

\newcommand{\Nat}{\mathbb{N}}
\newcommand{\ZZ}{\mathbb{Z}}
\newcommand{\CC}{\mathbb{C}}
\newcommand{\Hess}{\opname{Hess}}

\newcommand{\opname}[1]{\mbox{\rm #1}\,}

\newlength{\wex}  \newlength{\hex}
\def\ga{\alpha}     \def\gb{\beta}       \def\gg{\gamma}
       \def\gd{\delta}      
                         \def\vge{\varepsilon}
\def\gf{\phi}           
            \def\gl{\lambda}
\def\gm{\mu}        \def\gn{\nu}         
            
\def\gs{\sigma}

\def\Gg{\Gamma}           

\def\Gl{\Lambda}          
              
\newtheorem{Thm}{Theorem}[section]
\newtheorem{Lem}[Thm]{Lemma}

\newtheorem{Rems}[Thm]{Remarks}
\newtheorem{Rem}[Thm]{Remark}

\newcommand{\pf}{\noindent \mbox{{\bf Proof}: }}
\def\squarebox#1{\hbox to #1{\hfill\vbox to #1{\vfill}}}
\newcommand{\qed}{\hspace*{\fill}
      \vbox{\hrule\hbox{\vrule\squarebox{.667em}\vrule}\hrule}\smallskip}

\newcommand{\beq}{\begin{equation}}
\newcommand{\eeq}{\end{equation}}
\newcommand{\beqa}{\begin{eqnarray}}
\newcommand{\eeqa}{\end{eqnarray}}
\newcommand{\beqanl}{\begin{eqnarray*}}
\newcommand{\eeqanl}{\end{eqnarray*}}
\begin{document}
\renewcommand{\theequation}{\thesection.\arabic{equation}}
\newcommand{\mysection}[1]{\section{#1}\setcounter{equation}{0}}

\def\stackunder#1#2{\mathrel{\mathop{#2}\limits_{#1}}}
\newtheorem{theorem}{Theorem}
\newtheorem{lemma}[theorem]{Lemma}
\newtheorem{definition}[theorem]{Definition}
\newtheorem{corollary}[theorem]{Corollary}
\newtheorem{conjecture}[theorem]{Conjecture}
\newtheorem{remark}[theorem]{Remark}
\def\stackunder#1#2{\mathrel{\mathop{#2}\limits_{#1}}}

\author{Peter Kuchment \\
Mathematics and Statistics Department\\ Wichita State University\\
Wichita, KS 67260-0033, USA\\ kuchment@twsuvm.uc.twsu.edu
 \and Yehuda Pinchover \\ Department
of Mathematics\\ Technion - Israel Institute of Technology\\ Haifa
32000, Israel\\ pincho@techunix.technion.ac.il}
\title{Integral representations and Liouville theorems
for solutions of periodic elliptic equations}
\maketitle
\begin{abstract}
The paper contains integral representations for certain classes of
exponentially growing solutions of second order periodic elliptic
equations. These representations are the analogs of those
previously obtained by S.~Agmon, S.~Helgason, and other authors
for solutions of the Helmholtz equation. When one restricts the
class of solutions further, requiring their growth to be
polynomial, one arrives to Liouville type theorems, which describe
the structure and dimension of the spaces of such solutions. The
Liouville type theorems previously proved by M.~Avellaneda and
F.~-H.~Lin, and J.~Moser and M.~Struwe for periodic second order
elliptic equations in divergence form are significantly extended.
Relations of these theorems with the analytic structure of the
Fermi and Bloch surfaces are explained.\\[2mm]


\end{abstract}

\mysection{Introduction}\label{Sec1}

The topic of this paper stems from two sources. The
first of them are representation theorems for certain classes of
eigenfunctions of the Laplace operator in $\Real^n$, or
equivalently, of solutions of the Helmholtz equation
\begin{equation}
-\Delta u-k^2u=0 \quad \mbox{ in } \Real^n\, ,  \label{Helm}
\end{equation}
where $k\in \CC^{\, *}:=\CC\setminus\{0\}$. Such theorems for
arbitrary solutions of (\ref{Helm}) were obtained in $\Real^2$ and
in the hyperbolic plane by S.~Helgason \cite{H1,H2}, and in
$\Real^n$ by M.~Hashizume et al. \cite{HKMO}, M.~Morimoto
\cite{Mo}, and recently by S.~Agmon \cite{A2}. We remark that it
should also be possible to deduce similar results from the
L.~Ehrenpreis' {\em fundamental principle}. The zero set of the
symbol of the operator in the left hand side of (\ref{Helm}) is
\[
\Sigma =\left\{ \xi \in \CC^{\, n}|\,\xi ^2=k^2\right\} ,
\]
where $\xi^2=\sum_{j=1}^n \xi_j^2$. The L.~Ehrenpreis' ``fundamental
principle'' in the particular case of equation (\ref{Helm}) claims
that any solution of (\ref{Helm}) can be represented as a
combination (i.e., an integral with respect to the
parameter $\xi$) of the exponential solutions
\[
e_{\,\xi }(x):=\exp (i\xi \cdot x),\quad \xi \in \Sigma,
\]
where $\xi \cdot x=\sum_{j=1}^n \xi_jx_j$ (see the details and
more precise formulation in \cite{E} or \cite{Pa1}). The set
$\Sigma $ is an irreducible analytic subset of $\CC^{\, n}$, which
is uniquely determined when $k\neq 0$ by its spherical subset
\[
S=\left\{ \xi \in \CC^{\, n}|\,\xi =k\omega ,\,\omega \in
S^{n-1}\subset \mathbb{R}^n\right\} .
\]
Here $S^{n-1}$ denotes the unit sphere in $\Real^n$.
It is clear then that due to the uniqueness of analytic
continuation, the exponential representation of solutions $u(x)$
of (\ref{Helm}) should be reducible to one that involves the
solutions $e_{\,\xi }$ with $\xi \in S$ only. Namely, consider the
restriction mapping from functions analytic on the whole
characteristic variety $\Sigma$ to the sphere $S$. Due to the
irreducibility and the uniqueness of analytic continuation, this
mapping is one-to-one. Hence, there is a function space on the
sphere $S$ which is the isomorphic image of the space of all
analytic functions on $\Sigma$. It follows that any hyperfunction
(analytic functional) on $\Sigma$ can be rewritten as a functional
on $S$. Since the ``fundamental principle'' essentially expresses
all solutions of (\ref{Helm}) as applications of such analytic
functionals to the analytic family of exponential solutions, we
get our conclusion.

Now, depending on how fast the solution $u(x)$ grows at infinity,
the corresponding representing functional on $S$ is actually a
measure, a distribution, a hyperfunction, or a functional of a
more general kind.
For instance (see \cite{A2,HKMO,Mo}), an arbitrary solution in
$\Real^n$ can be represented as
\begin{equation}\label{rep1}
 u(x)=<\phi (\xi ),e_{\,\xi }(x)>,
\end{equation}
 where $\phi (\xi )$ is a functional on $S$ which belongs to the
dual space to the space $%
\mathcal{E}:=\lim_{R\to \infty }\mathcal{E}_R(S^{n-1})$. Here for
every $R>0$ the Hilbert space $\mathcal{E}_R(S^{n-1})$ is defined
as follows:
 \bea \mathcal{E}_R(S^{n-1}):=\left\{\right.\psi
\,\left|\right.\, && \hspace{-7mm} \psi (\omega )=\sum_{l,m}
a_{l,m}\frac{(R/2)^l}{\Gamma (l+\frac{n+1}{2})}Y_l^m(\omega ),
\mbox{ s.t. } \nonumber\\ &&\hspace{10mm} \|\psi
\|_{\mathcal{E}_R}:=(\sum_{l,m} |a_{l,m}|^2)^{1/2}<\infty
\left\}\right.,\nonumber \eea
 where $Y_l^m(\omega )$ denote
spherical harmonics, and $\mathcal{E}(S^{n-1})$ is equipped with
the inductive limit topology of $\lim_{R\to \infty
}\mathcal{E}_R(S^{n-1})$.

The representation (\ref{rep1}) can be formally rewritten as
\[
u(x)=\int_S e_{\,\xi }(x)d \phi (\xi ).
\]
where $\phi(\xi)$ is a suitable functional.

The functional $\phi $ is a hyperfunction (analytic functional) on
$S$ if and only if
for arbitrary $%
\varepsilon >0$ the solution $u(x)$ grows not faster than
\[
O(\exp ((|Im \, k|+\varepsilon )|x|))
\]
(see \cite{A2}). One can also describe other classes of solutions,
for instance, solutions which are represented by a distribution or
a measure (see \cite{A1a,A3,A2,LP} and the references therein).


As we have already mentioned, these results could be probably
extracted from the ``fundamental principle'' \cite{E, Pa1}. The
crucial factors are that $S$ is sufficiently massive and $\Sigma $
is irreducible, so $S$ determines $\Sigma $ uniquely. Besides, $S$
is a rather simple analytic manifold. These features allow more or
less explicit descriptions of the needed spaces of test functions
and functionals. It is easy to understand that if $\Sigma $ were
reducible, it would not be possible to obtain the representation
of {\em all} solutions using only $\xi \in S$. The reason is that
the solution $e_{\,\xi }$ with $\xi $ that belongs to a component
not touching $S$ would not be representable this way. On the other
hand, if one wants to deal only with solutions growing not faster
than $O(\exp (|Im \, k|+\varepsilon )|x|)$ for all $\varepsilon
>0$, then the irreducibility of $\Sigma $ is not needed. In this
case, it is only required that $\Sigma $ is irreducible in a
vicinity of $S$, so other components of $\Sigma $ do not meet $S$.

The ``fundamental principle'' was extended in \cite{Ku} to
solutions of certain growth (for instance, of exponential growth)
of elliptic and hypoelliptic {\em periodic} equations (see also
the extensions of the results of \cite{Ku} provided in
\cite{Pa2}). The role of the exponential solutions is played here
by the so called {\em Floquet-Bloch} solutions (see Definition
\ref{deffloq}), and an analog of the characteristic manifold
$\Sigma $ is the variety $F$ sometimes called the {\em Fermi
surface} (see Definition \ref{defFermi}). This raises the hope of
finding representations similar to the ones discussed above for
the more general case of a second order elliptic operator with
periodic coefficients. This is not straightforward, however, due
to several reasons. First of all, it is not that clear what should
be a natural analog of $S$. An appropriate variety, as we explain
later, is provided for second order equations by the analysis of
the cone of positive solutions done by S. Agmon and by V. Lin and
Y. Pinchover (see \cite{A1a,LP,Ku}, and the references therein).
The disadvantage is that the whole consideration must be done {\em
below} the spectrum of the operator (more precisely, below {\em
the generalized principal eigenvalue} $\Lambda_0$, see
(\ref{genev})). Secondly, proving the irreducibility of $F$
happens to be a very hard nut to crack (this problem arises also
in direct and inverse spectral problems, see for instance
\cite{BKT,GKT,KT,KV1,KV2}). Fortunately enough, by appropriately
restricting the growth of the solutions, one can sometimes work
near a single irreducible component, and hence avoid proving the
irreducibility of $F$. Consequently, we prove a representation
theorem (Theorem \ref{hyperfunction}) that characterizes all the
solutions which have integral expansion into {\em positive} Bloch
solutions with a hyperfunction as a ``measure".

The ``fundamental principle'' also suggests a point of view that
is crucial for establishing representation theorems for solutions
of equations with constant or periodic coefficients. Namely, it is
to one's advantage to treat solutions of the original equation in
the dual sense, i.e., as functionals on appropriate spaces that
are orthogonal to the range of the dual operator. We adopt this
approach throughout the paper. \vspace{3mm}

If one attempts now to further restrict the growth of solutions
and considers the problem of the structure of all polynomially
growing (or bounded) solutions, one arrives at the second topic of
our study, Liouville type theorems. The classical Liouville
theorem characterizes the space of all harmonic functions in
$\Real^n$ of polynomial growth of order $N$. The validity of an
analog of this classical theorem has been studied in many
situations (see for instance \cite{CM, L, LW} for recent results,
surveys, and further references). An interesting case was considered
by M.~Avellaneda and F.~-H.~Lin \cite{AL}, and also  by J.~Moser
and M.~Struwe \cite{MS}. In these papers the authors dealt with
polynomially growing solutions of a second order elliptic equation
$Lu=0$ in divergence form with periodic coefficients and obtained
a comprehensive answer (for related results see also \cite{CM,L}
and the references therein). More precisely, using the formalism of
homogenization theory \cite{BLP, JKO}, it was proved that any
solution $v$ of the equation $Lu=0$ in $\mathbb{R}^n$ of
polynomial growth is representable as a finite sum of the form
\begin{equation}\label{polyn}
v(x)=\sum\limits_{j=(j_1,\ldots,j_n)\in \ZZ_{+}^n}x^jp_j(x),
\end{equation}
where the functions $p_j(x)$ are periodic with respect to the
group of periods of the equation. Moreover, the space of all
solutions of the equation $Lu=0$ of polynomial growth of order at
most $N$ is of dimension $h_{n,N}$, where
\begin{equation}\label{hn}
h_{n,N}:=
\left(\begin{array}{c}n+N\\N\end{array}\right)-
\left(\begin{array}{c}n+N-2\\N-2\end{array}\right)
\end{equation}
is the dimension of the space of all harmonic polynomials of degree
$\leq N$ in $n$ variables.
We will also use the notation
\begin{equation}
q_{n,N}:=
\left(\begin{array}{c}n+N\\N\end{array}\right)
\label{qn}
\end{equation}
for the dimension of the space of all polynomials of degree at
most $N$ in $n$ variables. Notice that $q_{n-1,N}$ also coincides
with the dimension of the space of all homogeneous polynomials of
degree $N$ in $n$ variables, so in particular, $
h_{n,N}=q_{n-1,N-1}+q_{n-1,N}$.

We remark that the method of \cite{AL,MS} can be slightly modified
to provide an extension of this Liouville theorem for general
second order elliptic equations with periodic coefficients under
the assumption that the generalized principal eigenvalue is zero
(see Appendix \ref{Appendix2} and also the recent paper
\cite{LW1}, where a partial result of this type was
independently obtained).

One can make an observation that these Liouville theorems
are actually of the same nature as the
representation theorems discussed above. In this case the analog
of the set $S$ is the single point $\xi =0$ and the representing
functional $\phi $ is a distribution supported at this point.
\vspace{2mm}

Let us recall the following standard notion of Floquet theory (see
\cite{Ea, Ku, RS}):
\begin{definition}\label{deffloq}
{\em A solution $u(x)$ representable as a finite sum of the form
\begin{equation}\label{floquetsol}
u(x)=e^{ik\cdot x}\left( \sum\limits_{j=(j_1,\ldots ,j_n)\in \ZZ%
_{+}^n}x^jp_j(x)\right)
\end{equation}
with functions $p_j(x)$ periodic with respect to the group of
periods of the equation is called a {\em
Floquet solution with a quasimomentum $k\in \CC^n$}.
Here $k\cdot x=\sum k_l x_l$
The maximum value of
$|j|=\sum\limits_{l=1}^nj_l$ in the representation
(\ref{floquetsol}) is said to be the {\em order} of the Floquet
solution. Floquet solutions of zero order are called {\em Bloch
solutions}.}
\end{definition}

One sees that the representation (\ref{polyn}) corresponds to a
Floquet solution with a zero quasimomentum. A Liouville  theorem
of the type mentioned above implies in particular that the only
{\em real} quasimomentum that can occur for the equation under
consideration is $k=0$ (modulo the action of the lattice
reciprocal to the group of periods). We show in the present paper
that the finiteness of the set of {\em real} quasimomenta for a
periodic elliptic equation is equivalent to the finite
dimensionality of the spaces of solutions having a given
polynomial growth and to their representation similar to, albeit
more general than (\ref{polyn}). This statement is very general
and holds for any periodic elliptic equation (it is also true
for hypoelliptic equations and systems, although we not address
them here). If some additional
information is available on the analytic behaviour of the
dispersion relations, one can find the {\em exact} dimensions of
these spaces (see Theorem \ref{Liouville}). We present in Theorem
\ref{applicat} some classes of second order equations
(including Schr\"{o}dinger, magnetic Schr\"{o}dinger, and general
second order elliptic equations with real periodic coefficients)
for which one can achieve all these sharp results. We show that
the problem of calculating the dimensions of the spaces of Floquet
solutions of a given polynomial growth reduces to a purely
function theoretic question and can be resolved in a very general
setting (Theorem \ref{Fl_dimen}).

The proofs of the results of this paper are largely dependent upon
the techniques of the Floquet theory developed in \cite{Ku}.
\vspace{2mm}

The outline of the paper is as follows. The next section
introduces necessary notations and preliminary results from the
Floquet theory and the theory of positive solutions of periodic
elliptic equations. In particular, we obtain a new general result
(Theorem \ref{Fl_dimen}) on the dimensions of the spaces of
Floquet solutions, which plays crucial role in our approach to
Liouville theorems. Section \ref{Sec3} contains the proof of the
integral representation (Theorem \ref{hyperfunction}) analogous to
Theorem 5.1 in \cite{A2}. In Section \ref{Sec4}, we discuss
Liouville type theorems. In particular, Theorems \ref{Liouville}
and \ref{applicat} are established. In order to make the reading
of the paper easier, we postpone the proofs of all the technical
lemmas to Section \ref{Sec5}. Some conclusions and remarks are
provided in Section \ref{remarks}. The Appendix contains
an alternative proof of a part of Theorem \ref{applicat} using
the homogenization technique similar to the one used in
\cite{AL,MS}.

Results of this paper related to Liouville theorems were presented
in March 2000 at the University of Toronto and at the Weizmann
Institute.

When the paper was being prepared for submission, P.~Li informed the
authors that the statement of the
third part of Theorem \ref{applicat} for the special case
of an operator of the form
$L=-\sum a_{ij}(x)\partial _{i}\partial _{j}$
was simultaneously and independently obtained in \cite{LW1}
using homogenization formalism .

 \mysection{Notations and
preliminary results}\label{notations}

Consider a linear (scalar) elliptic partial differential operator
$P(x,D)$ of order $m$ in $\Real^{n}$, $n\geq 2$ (in some parts of
the paper we will restrict the class of operators further). Here
we employ the standard notation $D=\frac 1i\frac
\partial {\partial x}$. The ellipticity is understood in the sense
of the nonvanishing of the principal symbol $P_{m}(x,\xi)$ of the
operator $P$ for all $\xi \in \Real^{n}\setminus 0$. The dual
operator (the formal adjoint) $P^{*}$ has similar properties. Here
we use the duality provided by the bilinear (rather than the
sesquilinear) form
\[
<g,f>=\int\limits_{\Real^n}f(x)g(x)dx.
\]

We assume that the coefficients of $P$ are smooth and
periodic with respect to a lattice $\Gamma $ in $\Real^n$.
The smoothness condition can be significantly reduced (see
the Section \ref{remarks}). In fact, so far we only need
that both operators $P$ and $P^{*}$ define Fredholm
mappings between the Sobolev space $H^{m}$ and $L_2$ on the
torus $\mathbb{T}^n=\Real^{n} / \Gamma$.

An additional condition is required that would guarantee the
discreteness of the spectrum of the ``shifted" operators
$P(x,D+k)$ on the torus $\mathbb{T}^n$ for all $k\in
\mathrm{C}^n$. We need to exclude the possible pathological
situation when the spectrum of $P$ on the torus coincides with the
whole complex plane (like in the case of the operator $\exp(i\,
\phi)\, d/d \phi$ on the circle). For instance, self-adjointness
of $P$ could be such a condition. Another example is second order
uniformly elliptic operators of the form (\ref{operator}). For
more sufficient conditions see for example \cite{Ag}.

In what follows, the particular choice of the lattice is
irrelevant and can always be reduced to the case $\Gamma =\ZZ^n$,
which we will assume from now on. We will always use the word
``periodic" in the meaning of ``$\Gamma$-periodic".

We denote by $K=[0,1]^n\,$ the standard fundamental domain (the
{\it Wigner -Seitz cell}) of the lattice $\Gamma =\ZZ^n$, and by
$B=[-\pi ,\pi ]^n$ the {\it first Brillouin zone}, which is a
fundamental domain of the reciprocal (dual) lattice $\Gamma
^{*}=\left( 2\pi \ZZ\right) ^n$. We identify $\Gamma$-periodic
functions, in the natural way, with functions on ${\mathbb
T}^{n}$.

We introduce now the set that plays the role of the characteristic
variety $\Sigma $ discussed in the introduction.
\begin{definition} \label{defFermi}
{\em The (complex) {\em Fermi surface} $F_P$ of the operator $P$ (at
the zero energy level) consists of all vectors $k\in \CC^{\, n}$
(called {\em quasimomenta}) such that the equation $Pu=0$ has a
nonzero Bloch solution $u(x)=e^{ik\cdot x}p(x)$, where
 $p(x)$ is a $\Gamma $-periodic function.}
\end{definition}

It would be useful later on to realize that in this definition the
positivity of the solution is not required, and in fact the
solution is usually complex.

In many cases,  it is convenient to introduce a spectral parameter
$\lambda$. This leads to the notion of the {\it Bloch variety}:
\begin{definition} \label{defBloch}
{\em The (complex) {\em Bloch variety} $B_P$ of the operator $P$
consists of all pairs $(k,\lambda) \in \CC^{\, n+1}$
such that the equation $Pu=\lambda \, u$ has a
nonzero Bloch solution $u(x)=e^{ik\cdot x}p(x)$
with the quasimomentum $k$.}
\end{definition}

It is clear that the Fermi surface is just the projection onto the
$k$-space of the intersection of the Bloch variety with the
hyperspace $\lambda=0$.

One can consider the Bloch variety $B_P$ as the graph of a
(multivalued) function $\lambda(k)$, which is usually called the
{\it dispersion relation}. Then the Fermi surfaces become the
level surfaces of the dispersion relation. Since the spectra of
all operators $P(x,D+k)$ on the torus $\mathbb{T}^n$ are discrete,
we can single out continuous branches $\lambda_j$ of this
multivalued dispersion relation. These branches are usually called
the {\it band functions} (see \cite{RS, Ku}).

The following analyticity property of the Fermi and Bloch varieties
is important:
\begin{lemma}
\label{Fermi} \cite[Theorems 3.1.7 and 4.4.2]{Ku} The Fermi and
Bloch varieties are the sets of all zeros of entire functions of a
finite order in $\CC^{\, n}$ and $\CC^{\, n+1}$, respectively.
\end{lemma}

Another property of the Bloch and Floquet varieties that we will
need later is the relation between the corresponding varieties of
the operators $P$ and $P^*$.
\begin{lemma}
\label{dual} \cite[Theorem 3.1.5]{Ku} A quasimomentum $k$ belongs
to $F_{P^*}$ if and only if $-k \in F_{P}$. Analogously,
$(k,\lambda)\in B_{P^*}$ if and only if $(-k,\lambda)\in B_{P}$.
In other words, the dispersion relations $\lambda(k)$ and
$\lambda^*(k)$ for the operators $P$ and $P^*$ are related as
follows:
\begin{equation}
\lambda^*(k)=\lambda(-k).
\label{dual_disp}
\end{equation}
\end{lemma}

We note that the Fermi surface $F_P$ is periodic with respect to
the reciprocal lattice $\Gamma ^{*}=(2\pi \ZZ)^n$. Therefore, it
is sometimes useful to factor out the periodicity by considering
the (analytic) exponential mapping $\rho :\CC^{\, n}\rightarrow
(\CC^{\, *})^n$, where
\[
z=\rho (k)=\rho (k_1,\ldots,k_n)=(\exp ik_1,\ldots,\exp ik_n).
\]
This mapping can be identified in a natural sense with the
quotient map $\CC^{\, n}\rightarrow \CC^{\, n}/\Gamma ^{*}$.
We also introduce the complex torus
\begin{equation}
T=\rho (\Real^n)=\left\{ z\in \CC^{\, n}|\,\left| z_j\right|
=1,\,j=1,2,\dots,n\right\} .  \label{torus}
\end{equation}

\begin{definition}
\label{def_Floquet} {\em The image $\Phi_P=\rho (F_P)$ of the
Fermi surface $F_P$ under the mapping $\rho$ is called the {\em
Floquet surface} of the operator $P$.}
\end{definition}

The reader familiar with the Floquet theory immediately recognizes
the Floquet surface as the set of all Floquet multipliers of the
equation $Pu=0$.

The main tool in the Floquet theory is an analog of the Fourier
transform (see \cite[Section 2.2]{Ku}, \cite{RS}), which we will
call the {\it Floquet transform} ${\cal U}$ (it is sometimes also
called the {\it Gelfand transform}):
\begin{equation}
f(x)\rightarrow {\cal U}f(z,x)=\sum_{\gamma \in \Gamma }f(x-\gamma
)z^\gamma ,\quad z\in (\CC^{\, *})^n,  \label{Floquet}
\end{equation}
where we denote
 $z^\gamma =z_1^{\gamma _1}z_2^{\gamma_2}\ldots z_n^{\gamma _n}$ .

It is often convenient to use for the Floquet transform ${\cal U}$
the quasimomentum coordinate $k$ instead of the multiplier $z$, where
$$
z= \rho(k)=\left( \exp ik_1, \ldots , \exp ik_n \right).
$$

We need to recall now some definitions from \cite{Ku}. For a point
$z\in (\CC^{\, *})^n$, we denote by $E_{m,z}$ the closed subspace
of the Sobolev space $H^m(K)$ formed by the restrictions of
functions $v\in H_{loc}^m(\Real^n)$ that satisfy the Floquet
condition $v(x+\gamma )=z^\gamma v(x)$ for any $\gamma \in \Gamma
$. One can show (see Theorem 2.2.1 in \cite{Ku}) that
\begin{equation}
{\cal E}_m:=\mathrel{\mathop{\cup }\limits_{z\in (\CC^{\, *})^n}}%
E_{m,z}  \label{bundle}
\end{equation}
forms a holomorphic sub-bundle of the trivial bundle $(\CC^{\,
*})^n\times H^m(K)$. As any infinite dimensional analytic Hilbert
bundle over a Stein domain, it is trivializable (see Theorems
1.3.2, 1.3.3, and 1.5.23 in \cite{Ku}). One can also notice that
for $m=0$ the bundle ${\cal E}_0$ coincides with the whole
$(\CC^{\, *})^n\times L^2(K)$.

The following standard auxiliary result for the transform ${\cal
U}$ collects several statements from Theorem {\small{\sf XIII}}.97
in \cite{RS} and Theorem 2.2.2 in \cite{Ku}:
\begin{lemma}\label{Plancherel}
\begin{enumerate}
\item  For any nonnegative integer $m$ the operator
\[
{\cal U}:H^{m}(\Real^n)\rightarrow L^2(T,{\cal E}_m)
\]
is an isometric isomorphism, where $L^2(T,{\cal E}_m)$ denotes the
space of square integrable sections over the complex torus $T$ of
the bundle ${\cal E}_m$, equipped with the natural topology of a
Hilbert space.
\item Let the space
 $$ \Theta ^m=\{f\in H_{loc}^m(R^n)|\,
\underset{\gamma \in \Gamma }{sup}||f||_{H^m(K+\gamma )}\exp
(b|\gamma|)<\infty \,,\; \forall b>0\} $$ be equipped with the
natural Fr\'{e}chet topology. Then
\[
{\cal U}:\Theta^{m}\rightarrow \Gamma((\CC^{\, *})^{n},{\cal E}_m)
\]
is an isomorphism, where $\Gamma((\CC^{\, *})^{n},{\cal E}_m)$ is
the space of all analytic sections over $(\CC^{\, *})^{n}$ of the
bundle ${\cal E}_m$, equipped with the topology of uniform
convergence on compacta.
\item Let the elliptic operator $P$ be of order $m$. Then under
the transform ${\cal U}$ the operator
$$
P: H^{m}(\Real^n) \rightarrow L^{2}(\Real^n)
$$
becomes the operator of multiplication by a holomorphic
Fredholm morphism $P(z)$ between the fiber bundles
${\cal E}_m$ and ${\cal E}_0$.
Here $P(z)$ acts on the fiber of ${\cal E}_m$ over
the point $z \in T$ as the restriction to this fiber of
the operator $P$ acting between $H^m(K)$ and $L^2(K)$.
\end{enumerate}
\end{lemma}

Here is another standard way of looking at the morphism $P(z)$.
Let $z=\exp ik$, then commuting with the exponent $\exp ik \cdot
x$ one can (locally) trivialize the bundle ${\cal E}_m$ reducing
it to the trivial bundle with the fiber $H^{m}({\mathbb T}^n)$,
where as before ${\mathbb T}^n=\Real^n / \Gamma$. At the same time
the operator $P(z)$ takes the form $P(x,D+k)$ between Sobolev
spaces on the torus ${\mathbb T}^n$.

Let us discuss the structure of the Floquet solutions (see
Definition \ref{deffloq}) and of functions of Floquet type
(\ref{floquetsol}) in general. For illustration, consider the
constant coefficient case, where the role of the Floquet solutions
is played by the exponential polynomials $$ e^{ik\cdot
x}\sum\limits_{\left| j\right| \leq N}p_jx^j. $$ It is well known
that, considered as distributions, all such functions are Fourier
transformed into distributions supported at the point
 $\left(-k\right)$. Moreover, the converse statement is also true. A simple but
extremely important and relatively unnoticed observation is that
under the Floquet transform, each Floquet type function of the
form (\ref{floquetsol}) corresponds, in a similar way, to a
(vector valued) distribution supported at the quasimomentum
$\left( -k\right) $. We collect below this fact and some other
previously known properties of Floquet solutions, as well as a new
result on the dimensions of the spaces of such solutions, which
will play the crucial role in establishing the Liouville type
theorems.

First of all, every Floquet type function $u$ (see
(\ref{floquetsol})), being of exponential growth, determines a
(continuous linear) functional on the space $\Theta ^0 $. If,
additionally, it satisfies the equation $Pu=0$ for a periodic
elliptic operator of order $m$, then as such a functional it is
clearly orthogonal to the range of the dual operator $P^{*}:\Theta
^m\rightarrow \Theta ^0$. According to Lemma \ref{Plancherel},
after the Floquet transform any such functional becomes a
functional on $\Gamma \left( \left( \CC^{\, *}\right) ^n,{\cal
E}_0\right)$, which is orthogonal to the range of the Fredholm
morphism $P^{*}(z):{\cal E}_m\rightarrow {\cal E}_0$ generated by
the dual operator $P^{*}:\Theta ^m\rightarrow \Theta ^0$. We are
now ready to formulate the following auxiliary result.
\begin{lemma}\label{Floq_struct}
\begin{enumerate}
\item
A continuous linear functional $u$ on $\Theta ^0$ is generated by
a function of the Floquet form (\ref{floquetsol}) with a
quasimomentum $k$ if and only if after the Floquet transform it
corresponds to a functional on $\Gamma \left( \left( \CC^{\,
*}\right) ^n,{\cal E}_0\right) $ which is a distribution $\phi $
that is supported at the point $\nu =\exp (-ik)$, i.e. has the
form $$ \left\langle \phi ,f\right\rangle =\sum_{\left| j\right|
\leq N}\left\langle q_j,\left. \frac{\partial ^{\left| j\right|
}f}{\partial z^j}\right| _{\,\nu }\right\rangle , \,f\in \Gamma
\left( \left( \CC^{\, *}\right) ^n,{\cal E}_0\right), $$ where
$q_j\in L^2(K)$. The orders $N$ of the Floquet function
(\ref{floquetsol}) and of the corresponding distribution
$\phi $ are the same.
\item
Let $a_k$ be the dimension of the kernel of the operator $$
P(x,D+k):H^m({\mathbb T}^{n})\rightarrow L^2({\mathbb T}^{n}). $$
Then the dimension of the space of Floquet solutions of the
equation $Pu=0$ of order at most $N$ with a quasimomentum $k$ is
finite and does not exceed $a_{k}q _{n,N}$.
\end{enumerate}
\end{lemma}

The estimate on the dimension given in the second part of Lemma
\ref{Floq_struct} is very crude and in many cases can be
significantly improved. In fact, as the following theorem shows,
we obtain an explicit formula for the dimension of the space of
Floquet solutions with a given quasimomentum in the case of a
simple eigenvalue. This theorem seems to be new and constitutes
the crucial part of the Liouville theorem proved in Section
\ref{Sec4} (Theorem \ref{Liouville}).

In order to formulate this result, we need to prepare some notions
and notations.
\begin{definition}\label{Qharm} {\em Let $Q$ be a homogeneous polynomial
in $n$ complex variables. A polynomial $p(x)$ in $\Real^n$ is
called {\em $Q$-harmonic}, if it satisfies the differential
equation $Q(D)p=0$.}
\end{definition}

Let ${\cal P}$ denote the vector space of all polynomials in $n$
variables, and let $P_l$ be the subspace of all homogeneous
polynomials of degree $l$. Denote by ${\cal
P}_N=\bigoplus\limits_{l=0}^NP_l$ the subspace of all polynomials
of degree at most $N$. So, ${\cal P}=\bigoplus\limits_{l=0}^\infty
P_l$. If $Q(k)$ is a nonzero homogeneous polynomial of degree $s$,
then the differential operator $Q(D):P_{l+s}\rightarrow P_l$ is
surjective for any $l$ (this simple statement will also follow
from the proof of the theorem below). Hence, the mapping
$Q(D):\cal{P}\rightarrow \cal{P}$ has a (nonuniquely defined)
linear right inverse $R$ that preserves the homogeneity of
polynomials.

\begin{theorem}
\label{Fl_dimen} Assume that zero is an eigenvalue of algebraic
multiplicity $1$ of the operator
$P(x,D+k_{0}):H^m(\mathbb{T}^n)\rightarrow L^2(\mathbb{T}^n)$ on
the torus $\mathbb{T}^n$. Let $\lambda (k)$ be an analytic
function in a neighborhood of $k_0$ such that $\lambda (k)$ is a
simple eigenvalue of the operator
$P(x,D+k):H^m(\mathbb{T}^n)\rightarrow L^2(\mathbb{T}^n)$ and
$\lambda(k_{0})=0$. Let $$ \lambda (k)=\sum\limits_{l=l_0}^\infty
\lambda _l(k-k_0) $$ be the Taylor expansion of $\lambda (k)$
around the point $k_0$ into homogeneous polynomials such that
$\lambda _{l_0}$ is the first nonzero term of this expansion. Then
for any $N\in \Nat$ the dimension of the space of Floquet
solutions of the equation $Pu=0$ in $\Real^n$ of order at most $N$
and with the quasimomentum $k_0$
is equal to the dimension of the space of all $\lambda
_{l_0}$-harmonic polynomials of degree of at most $N$. Moreover,
given a linear right inverse $R$ of the mapping $\lambda
_{l_0}(D):\cal{P}\rightarrow \cal{P}$ that preserves homogeneity,
one can construct an explicit isomorphism between these spaces.
\end{theorem}
{\bf Proof}. It is sufficient to consider the case $k_0=0$, since
the general case reduces to this by a change of variables.
Consider the operator family $$A(k)=P^{*}(x,D-k)-\lambda
(-k):H^m(\mathbb{T}^n)\rightarrow L^2(\mathbb{T}^n) $$ which is
analytic in a neighborhood of $0$. At each point $k$ of this
neighborhood $A(k)$ has by the construction a one-dimensional
kernel. Then, according to Theorem 1.6.13 in \cite{Ku}, there
exists an analytic non-vanishing vector $\psi (k)\in KerA(k)$. In
other words, $P^{*}(x,D-k)\psi (x,k)=\lambda (-k)\psi (x,k)$. Let
us choose a closed complementary subspace $M$ to $KerA(0)$ in
$H^m(\mathbb{T}^n)$. Then it is complementary to $Ker\,A(k)$ in a
neighborhood of $0$. Since $P^{*}(x,D)$ has zero kernel on $M$ and
is Fredholm, we conclude that $P^{*}(x,D-k)$ has zero kernel on
$M$ for all $k$ in a neighborhood of $0$. We denote by $\Pi (k)$
the closed subspace in $L^2(\mathbb{T}^n)$ defined as $\Pi
(k)=P^{*}(x,D-k)(M)$. Applying Theorem 1.6.13 of \cite{Ku} again,
we conclude that $\Pi (k)$ depends holomorphically on $k$ in a
neighborhood of $0$ (i.e., forms a Banach bundle) and hence it is
complementary to $Ker\,A(k)$. Representing now the operator
$P^{*}(x,D-k)$ in the block form according to the decompositions
$$ H^m(\mathbb{T}^n)=M\oplus Ker\,A(k) $$ and $$
L^2(\mathbb{T}^n)=\Pi (k)\oplus Ker\,A(k), $$ we get $$
P^{*}(x,D-k)= \left(
\begin{array}{cc}
B(k)&0 \\
0&\lambda (-k)
\end{array}
\right), $$ where $B(k)$ is an analytic invertible
operator-function between $M$ and $\Pi (k)$. If we now have a
functional $\phi $ on $\Gamma \left( \left( \CC^{\, *}\right)
^n,{\cal E}_0\right) $ supported at $v=\exp (0)$, such that it is
orthogonal to the range of the operator of multiplication by
$P^{*}(k)$, then it must be equal to zero on all sections of the
bundle $\Pi (k)$. This means that the restriction of such
functionals to the sections of the one-dimensional bundle
$Ker\,A(k)$ is an one-to-one mapping. This reduces the problem to
the following scalar one: find the dimension of the space of all
distributions of order $N$ supported at the origin such that they
are orthogonal to the ideal generated by $\lambda (-k)$ in the
ring of germs of analytic functions. One can change variables to
eliminate the minus sign in front of $k$. Due to the finiteness of
the order of the distribution, the problem further reduces to the
following: find the dimension of the cokernel of the mapping
$$\Lambda_N
:{\cal P}_N\rightarrow {\cal P}_N\,,$$
where $\Lambda_N(p)$ for $p\in {\cal P}_N$ is the Taylor polynomial
of order $N$ at $0$ of the product $\lambda (k)p(k)$. Let us write
the block matrix $\Lambda_{ij}$ of the operator $\Lambda_N $ that
corresponds to the decomposition ${\cal
P}_N=\bigoplus\limits_{l=0}^NP_l$. It is obvious that $\Lambda
_{ij}=0$ for $i-j<l_0$. For $i-j\geq l_0$ the entry $\Lambda
_{ij}$ is the operator of multiplication by $\lambda _{i-j}$
acting from $P_j$ into $P_i$. Since $\lambda_{l_0} \neq 0$, for
$i-j= l_0$ the operator $\Lambda _{ij}$ of multiplication by
$\lambda_{l_0}$ has zero kernel. Being interested in the cokernel
of $\Lambda_N $, we need to find the kernel of the adjoint matrix
$\Lambda_N ^{*}$. The adjoint matrix acts in the space
$\bigoplus\limits_{l=0}^NP_l^{*}$, where $P_l^{*}$ can be
naturally identified with the space of linear combinations of the
derivatives of order $l$ of the Dirac's delta-function at the
origin. Here we have $\Lambda _{ij}^{*}=0$ for $j-i<l_0$, and for
$j-i\geq l_0$ the entry $\Lambda _{ij}^{*}$ is the dual to the
operator of multiplication by $\lambda _{j-i}$ acting from $P_i$
into $P_j$. In particular, since for $j-i= l_0$ the latter
operator is injective, we conclude that the operators $\Lambda
_{ij}^{*}$ are surjective. This enables one to find the dimension
of the kernel of the matrix $\Lambda_N ^{*}$ and even to describe
its structure. Namely, let $$ \psi =\left( \psi _0,...,\psi
_N\right) \in \bigoplus\limits_{l=0}^NP_l^{*} $$ be such that
$\Lambda ^{*}\psi =0$. Due to the triangular structure of
$\Lambda_N ^{*}$, it is easy to solve this system. Indeed, it can
be written as follows: $$ \sum\limits_{j\geq i+l_0}\Lambda
_{ij}^{*}\psi _j=0,\,i=0,...,N-l_0. $$ Taking the Fourier
transform, we can rewrite this system in the form $$
\sum\limits_{j\geq i+l_0}\lambda _{j-i}(D)\widehat{\psi
_j}=0,\,i=0,...,N-l_0, $$ where $\widehat{\psi}$ denotes the
Fourier transform of $\psi$. Therefore, $\widehat{\psi _j}$ is a
homogeneous polynomial of degree $j$ in $\Real^n$. For $i=N-l_0$
we have $$\lambda _{l_0}(D)\widehat{\psi _N}=0. $$ This equality
means that $\widehat{\psi _N}$ can be chosen as an arbitrary
$\lambda _{l_0}$-harmonic homogeneous polynomial of order $N$.
Moving to the previous equation, we analogously obtain $$ \lambda
_{l_0}(D)\widehat{\psi _{N-1}}+\lambda _{l_0+1}(D)\widehat{\psi
_N}=0, $$ or $$
 \lambda _{l_0}(D)\widehat{\psi _{N-1}}=-\lambda _{l_0+1}(D)\widehat{\psi _N}.
$$ The right hand side is already determined, and the
nonhomogeneous equation, as we concluded before, always has a
solution, for instance $$
 -R\left( \lambda _{l_0+1}(D)\widehat{\psi _N}\right) .
$$
This means that
$$
\widehat{\psi _{N-1}}+R\left( \lambda _{l_0+1}(D)\widehat{\psi _N}\right)
$$
 is a $\lambda _{l_0}$-harmonic homogeneous polynomial of order $N-1$. We see that
the solution $\widehat{\psi _{N-1}}$ exists and is determined up
to an addition of any homogeneous $\lambda _{l_0}$-harmonic
polynomial of degree $N-1$. Continuing this process until we reach
$\widehat{\psi _{0}}$, we conclude that the mapping $$ \psi
=\left( \psi _0,...,\psi _N\right) \rightarrow \phi =\left( \phi
_0,...,\phi _N\right), $$ where $$ \phi _j=\widehat{\psi
_j}+R\sum\limits_{i>j}\lambda _{i-j+l_0}(D)\widehat{\psi _i} $$
establishes an isomorphism between the cokernel of the mapping
$\Lambda_N $ and the space of $\lambda _{l_0}$-harmonic
polynomials of degree at most $N$. This proves the theorem. \qed

In the cases of the simplest structures of the Taylor series, the theorem implies
the following:
\begin{corollary}
\label{specific} Under the hypotheses of Theorem \ref{Fl_dimen}
the following hold:
\begin{enumerate}
\item
If $k_0$ is a noncritical point of the band function $\lambda
(k)$, then the dimension of the space of Floquet solutions of the
equation $Pu=0$ in $\Real^n$ of order at most $N$ with a
quasimomentum $k_0$ is equal to the dimension $q_{n-1,N}$ of the
space of all polynomials of degree at most $N$ in $\Real^{n-1}$.
\item
If the Taylor expansion of the band function $\lambda (k)$ at a
point $k_0$ starts with a nondegenerate quadratic form, then the
dimension of the space of Floquet solutions of the equation $Pu=0$
in $\Real^n$ of order at most $N$ with a quasimomentum $k_0$ is
equal to the dimension $h_{n,N}$ of the space of harmonic (in the
standard sense) polynomials of degree at most $N$ in $\Real^n$. In
particular, this condition is satisfied at nondegenerate extrema.
\end{enumerate}
In both cases an isomorphism can be provided explicitly as in the previous theorem.
\end{corollary}
\pf 1. By our assumptions, the Taylor expansion of $\lambda (k)$
starts with a nonzero linear term $\lambda_{1}(k)=\sum_{j=1}^n
a_{j}k_{j}, \, a_{j} \in \CC$. The corresponding differential
operator is $$ \lambda _1(D)=-i\sum_{j=1}^n a_j\frac \partial
{\partial x_j}= \sum_{j=1}^n \alpha _j\frac
\partial {\partial x_j} +i\sum_{j=1}^n \beta _j\frac \partial {\partial
x_j}, $$ where $\alpha _j$ and $\beta _j$ are real. Consider first
the case when the vectors $\alpha =(\alpha _j)$ and $\beta =(\beta
_j)$ are collinear. Then $\lambda _1(D)$ becomes $\gamma _0\sum
\gamma _j\frac \partial {\partial x_j}$, where $\gamma _0\neq 0$
is a complex number and $\gamma =(\gamma _j)$ is a nonzero real
vector. A linear change of coordinate system brings $\lambda
_1(D)$ to the operator $\frac \partial {\partial x_1}$ (up to an
irrelevant constant factor). Thus, the $\lambda _1$-harmonic
polynomials are exactly those independent on $x_1$. Invoking
Theorem \ref{Fl_dimen}, we get our conclusion in this case.
Consider now the situation when $\alpha $ and $\beta $ are
linearly independent. Then a linear change of variables brings
$\lambda _1$ to the form $\partial /{\partial \overline{z}}$,
where $z=x_1+ix_2$. Since any polynomial in variables
$(x_1,...,x_n)$ is a polynomial of the same degree in
$(z,\overline{z},x_3,...,x_n)$, the $\lambda _1$-harmonic
polynomials are the ones depending on $(z,x_3,...,x_n)$ only (i.e.
the ones analytic in $z$). This again reduces the number of
variables to $n-1$. \vspace{2mm}

\noindent 2.  By our assumptions, the first nonzero homogeneous
term is a nondegenerate quadratic form $\lambda _2(k-k_0)$, which
is reducible to the sum of squares of coordinates by a linear
change of variables. Therefore, in the new coordinates $\lambda
_2(D)=-\Delta$. Using Theorem \ref{Fl_dimen}, we obtain the
desired result. \qed

In the remaining part of this section we restrict further the form
of the operator. Namely, we consider now second order operators
with {\it real} periodic coefficients of the form
\begin{equation}
L=-\sum_{i,j=1}^n a_{ij}(x)\partial _{i}\partial _{j}+\sum_{i=1}^n
b_i(x)\partial _i+c(x). \label{operator}
\end{equation}
It is assumed that the uniform ellipticity condition
\[
\sum_{i,j=1}^n a_{ij}(x)\zeta _i\zeta _j\geq
a\sum_{i=1}^n\zeta_i^2
\]
is satisfied for all $x, \zeta \in \Real^n$, where $a$ is a
positive constant.

For such operators, we introduce the function that will play the
crucial role in our considerations. Its properties were studied in
detail in \cite{A1a}, \cite {LP}, and \cite{Pr}. Consider the
function $\Lambda (\xi ):\Real^n\rightarrow \Real$ defined by the
condition that the equation
\[
Lu=\Lambda (\xi )u
\]
has a {\em positive} Bloch solution of the form
\begin{equation}\label{positiveBloch}
u_{\,\xi }(x)=e^{\xi \cdot x}p_{\,\xi }(x),
\end{equation}
where $p_{\,\xi }(x)$ is $\Gamma $-periodic.

\begin{lemma}
\label{Lambda-lemma}
\begin{enumerate}
\item  The value $\Lambda (\xi )$ is uniquely determined for any $\xi \in
\Real^n$.

\item  The function $\Lambda (\xi )$ is bounded from above, strictly
concave, analytic, and has a nonzero gradient at all points except
at its maximum point.

\item  Consider the operator
$$L(\xi )=e^{-\xi\cdot x}Le^{\xi\cdot x}=L(x,D-i\xi) $$ on the
torus $\mathbb{T}^n$.  Then $\Lambda (\xi )$ is the principal
eigenvalue of $L(\xi )$ with  a positive eigenfunction $p_{\,\xi
}$. Moreover, $\Lambda (\xi )$ is algebraically simple.

\item The Hessian of $\Lambda (\xi )$ is nondegenerate at all points.

\end{enumerate}
\end{lemma}

One should note that since the function $\Lambda(\xi)$ is analytic,
it is actually defined in a neighborhood of $\Real^n$ in $\CC^n$.
This remark will be used in what follows.

Let us denote
\begin{equation}
\Lambda_{0} =\max_{\xi \in \Real^n}\Lambda (\xi ).  \label{Lambda}
\end{equation}
It follows from \cite{A1a,LP} that an alternative definition of $\Lambda_{0}$
is
\begin{equation}\label{genev}
\Lambda_0= \sup\{\gl \in \Real \; |\; \exists u>0 \mbox{ such that }
(L-\gl)u= 0 \mbox{ in } \Real^n\},
\end{equation}
and that in the self-adjoint case $\Lambda_{0}$ coincides with
the bottom of the spectrum of the operator $L$. The common name
for $\Lambda_{0}$ is {\em the generalized principal eigenvalue} of
the operator $L$ in $\Real^n$.

We will often need to assume that $\Lambda_{0}$ is either
nonnegative or strictly positive. In the self-adjoint case such an
assumption  has a clear spectral interpretation. In the next
lemma, we provide some known sufficient conditions for the
nonnegativity or positivity of $\Lambda_{0}$ for operators of the
form (\ref{operator}).

\begin{lemma}
\label{Lambdazero}
Consider an operator $L$ of the form (\ref{operator})
\begin{enumerate}

\item $\Lambda _{0}\geq 0$ if and only if the operator $L$ admits a positive
(super)solution. This condition is satisfied in particular when
$c(x)\geq 0$.

\item $\Lambda _{0}\geq 0$ if and only if the operator $L$ admits a positive
solution of the form (\ref{positiveBloch}).

\item $\Lambda _{0}=0$ if and only if the equation $Lu=0$ admits
exactly one normalized positive solution in $\Real^n$.

\item If $c(x)=0$, then $\Lambda_{0} =0$ if and only if
$\int\limits_{{\mathbb T}^n}b(x)\psi (x)\,dx=0$, where $\psi $ is
the principal eigenfunction of $L^{*}$ on ${\mathbb T}^n$. In
particular, divergence form operators satisfy this condition.

\item Let $\xi\in \Real^n$, and assume that  $u_\xi(x)=e^{\xi \cdot x}p_\xi (x)$
and $u^*_{-\xi}$ are positive Bloch solutions of the equations
$Lu=0$ and $L^*u=0$, respectively. Denote by $\psi$ the periodic
function $u_\xi u^*_{-\xi}\,$. Consider the function
\[
\tilde{b}_i(x)= b_i(x)-2\sum_{j=1}^n a_{ij}(x)\{\xi_j
+(p_\xi(x))^{-1}\partial_{j}p_\xi(x)\},
\]
and denote
\[
\gamma =(\gamma _1,\ldots ,\gamma
_n):=(\int\limits_{\mathbb{T}^n}\tilde{b}_1(x)\psi (x)\,dx,\ldots
,\int\limits_{\mathbb{T}^n}\tilde{b}_n(x)\psi (x)\,dx).
\]
Then $\Lambda_{0} =0$ if and only if $\gamma =0$.

\end{enumerate}
\end{lemma}
Let us discuss also  some additional properties that will play an
important role in the sequel. Assume that $\Lambda_{0}\geq 0$.
Then Lemma \ref{Lambda-lemma} implies that the zero level set
\begin{equation}
\Xi =\left\{ \xi \in \Real^n|\;\Lambda (\xi )=0\right\}
\label{ksi}
\end{equation}
is either a strictly convex compact analytic surface in $\Real^n$
of dimension $n-1$ (this is the case if and only if $\Lambda_{0}>
0$), or a singleton  (this is the case if and only if
$\Lambda_{0}=0$). The manifold $\Xi$ consists of all $\xi \in
\Real^n$ such that the equation $Lu=0$ admits a positive Bloch
solution $u_{\,\xi}(x)=e^{\xi\cdot x}p_{\,\xi }(x)$. Moreover, the
set of all such positive Bloch solutions is the set of all {\em
minimal positive solutions} of the equation $Lu=0$ in
$\mathbb{R}^n$ \cite{A1a,LP}. It is also established that a
function $u$ is a positive solution of the equation $Lu=0$ in
$\mathbb{R}^n$ if and only if there exists a positive finite
measure $\gm$ on $\Xi$ such that
\[
u(x)=\int_\Xi u_{\,\xi }(x)d\gm(\xi ).
\]


We denote by $G$ the convex hull  of $\Xi$, and by
$\stackrel{\circ}{G}$ its interior. Note that if $\Lambda_{0}\geq
0$ then $\Lambda_{0}=0$ if and only if $\Xi=G$ and hence
$\stackrel{\circ}{G}=\emptyset $.

\begin{lemma}
\label{analytic} Suppose that $\Lambda_{0}>0$. There exists a
neighborhood $W$ of $G$ in $\CC^n$ and an analytic function
$$
W\ni\xi \mapsto p_{\xi}(\cdot) \in H^2(\mathbb{T}^n )
$$
such that for any $\xi \in W$ the function of $x$
\[
u_\xi (x)=\exp(\xi \cdot x)p_{\xi}(x)
\]
is a nonzero Bloch solution of the equation $Lu=\Lambda (\xi)u$
with a quasimomentum $-i \xi$. Moreover, one can choose the
function $p$ in such a way that it is positive for all $\xi\in
\Xi$.
\end{lemma}

Comparing the definitions of $\Xi $ and of the Fermi surface
$F_L$, it follows that
\[
-i\Xi \subset F_L.
\]
The next lemma specifies further the relation between these two varieties:
\begin{lemma}
\label{Fermi-lemma} Let $\Lambda_{0} \geq0$. Then
\begin{enumerate}
\item  The intersection of the complex Fermi surface $F_L$ with the tube
\begin{equation}
{\cal T}=\left\{ k\in \CC^{\, n}|\;Im \, k=(Im \, k_1,\dots ,Im \,
k_n)\in -G\right\} \label{tube}
\end{equation}
coincides with the union of the surface $-i\Xi $ with its
translations by the vectors of the reciprocal lattice $\Gamma
^{*}$, i.e. consists of vectors $k=-i\xi +\gamma $ where $\xi \in
\Xi $ and $\gamma \in \Gamma ^{*}$. Moreover, up to a
multiplicative constant, any nonzero Bloch solution with a quasimomentum
in the above intersection is a positive Bloch solution.

\item  If $\Lambda_{0} >0$, then the intersection of $F_L$ with a sufficiently
small neighborhood of $-i\Xi $ is a (smooth) analytic manifold
that coincides with the set of zeros of the function
$\Lambda(ik)$.
\end{enumerate}
\end{lemma}

Analogously to the definition of the Floquet surface
$\Phi=\Phi_L$, we define the surface
\begin{equation}
\Psi =\rho (-i\Xi )=\left\{ z\,|\; z=(\exp \xi_1,\ldots, \exp
\xi_n),\;\;\xi\in \Xi\right\}, \label{Psi}
\end{equation}
and the tubular domain
\begin{equation}\label{tubulard}
V=\rho({\cal T}),
\end{equation}
where ${\cal T}$ was defined in (\ref{tube}). The results of
Lemmas \ref{analytic} and \ref{Fermi-lemma} can be rephrased in
terms of these objects:

\begin{lemma}
\label{Lambda-lemma2} Let $\Lambda_0\geq 0$. Then
\begin{enumerate}
\item  $\Phi \cap V=\Psi $.

If $\Lambda_{0} >0$, then

\item  The intersection of $\Phi $ with a sufficiently small neighborhood of
$\Psi $ is a (smooth) connected analytic manifold.

\item  The intersections of $\Phi $ with neighborhoods of the tube $V$ form
a basis of neighborhoods of $\Psi $ in $\Phi $.

\item  For a sufficiently small neighborhood $\Phi _{\,\varepsilon }$ of
$\Psi $ in $\Phi $ there exists an analytic function
$p:\Phi_{\,\varepsilon }\rightarrow H^2(\mathbb{T}^n)$ such that
for any $z\in \Phi _{\,\varepsilon }$ the function of $x$
\[
u_z(x)=z^xp(z,x)
\]
is a nonzero Bloch solution of the equation $Lu=0$.
\end{enumerate}
\end{lemma}

We will also employ the following lemma:

\begin{lemma}
\label{constants}
Consider an operator $L$ of the form (\ref{operator})
\begin{enumerate}

\item Assume that $c(x)\geq 0$. Then the only solutions of the
equation $Lu=0$ of the type $\exp (ik\cdot x)p(x)$, where $k\in
\Real^n$ and $p$ is a $\Gamma$-periodic function are the
constants. If such a nontrivial solution exists, then $c(x)=0$,
and $\Gl(0)=0$ (i.e. $0\in \Xi$).

\item Suppose that the operator $L$ admits a positive periodic
supersolution $\psi \in C^{2,\ga}(\Real ^n)$. Assume that
$v(x)=\exp (ik\cdot x)p(x)$ is a nontrivial solution of the
equation $Lu=0$, where  $k\in \Real^n$ and $p$ is a periodic
function. Then there exists $C\in \CC$ such that $v=C\psi$, the
function $\psi$ is a positive periodic solution, and $\Gl(0)=0$.
\end{enumerate}
\end{lemma}

\mysection{Representation of solutions by
hyperfunctions}\label{Sec3}

The main result of this section (Theorem \ref{hyperfunction}) is
analogous to Theorem 5.1 in \cite{A2}, which characterizes the
class of solutions of the Helmholtz equation that can be
represented by means of hyperfunctions on $S$ (see also the
introduction of our paper).

In order to state it, we need to introduce a new object. Let us
denote by $h(\omega )$,  $\omega \in S^{n-1}$ the indicator
function of the convex set $G$. Namely,
\begin{equation}
h(\omega )=\mathrel{\mathop{\sup }\limits_{\xi\in G}}(\omega \cdot
\xi), \label{indicator}
\end{equation}
where $\omega \cdot \xi=\sum_{j=1}^n \omega _j\xi_j$ is the inner
product in $\Real^n$. The next Theorem will be stated in terms of
this function.

\begin{theorem} \label{hyperfunction}
Suppose that $\Lambda_{0} >0$. Let  $u$ be a solution of the
equation $Lu=0$ in $\Real^n$ satisfying for any $\varepsilon >0$
the estimate
\begin{equation}\label{estimate}
\left| u(x)\right| \leq C_{\,\varepsilon }\exp \left( \left(
h(x/\left| x\right| )+\varepsilon \right) \left| x\right| \right),
\end{equation}
where $C_{\,\varepsilon }$ is a constant depending only on
$\varepsilon$ and $u$. Then $u$ can be represented as
\begin{equation}
u(x)=<\mu (\xi ),u_{\,\xi }(x)>,  \label{represent}
\end{equation}
where $u_{\,\xi} $ is the analytic positive Bloch solution
corresponding to $\xi \in \Xi$ (see Lemma \ref{analytic}), and $
\mu (\xi )$ is a hyperfunction (analytic functional) on $\Xi $.
The converse statement is also true: for any hyperfunction $\mu $
on $\Xi $, the function $u(x)$ in (\ref{represent}) is a solution
of the equation $Lu=0$ in $\Real^n$ which satisfies the growth
condition (\ref{estimate}).
\end{theorem}
\begin{remark}
\label{pointwiseL2} {\em Using a standard elliptic argument it
follows that the pointwise growth condition (\ref{estimate}) is
equivalent to the growth condition
\begin{equation}\label{L2estimate}
u(x)\exp \left( -\left( h(x/\left| x\right| )+\varepsilon \right)
\left| x\right| \right)\in L^2(\mathbb{R}^n)\,.
\end{equation}
}\end{remark}

\noindent \mbox{{\bf Proof of Theorem \ref{hyperfunction}}:
}Assume first that a solution $u$ has the representation (\ref
{represent}). We need to prove that $u$ satisfies the growth
condition (\ref{estimate}). Due to the real analyticity of
$u_{\,\xi }$ with respect to $\xi $ and according to lemmas
\ref{analytic} and \ref{Fermi-lemma}, $u_{\,\xi }$ can be extended
to an analytic vector
function $u_{\,\xi }(x)=\exp (\xi \cdot x)p_{\xi}(x)$ on an $\varepsilon $%
-neighborhood $U_\varepsilon $ of $\Xi $ in $iF_{L}$. Since $\mu$ is a
hyperfunction (analytic functional) on $\Xi $, we have an estimate
\[
\left| u(x)\right| \leq C_{\,\varepsilon }\max_{\xi \in
U_{\,\varepsilon }}\left| u_{\,\xi }(x)\right| .
\]
Hence we have
\[
\left| u(x)\right| \leq C_{\,\varepsilon }\max_{\xi \in
U_{\,\varepsilon }}\left| e^{\xi \cdot x}\right| =C_{\,\varepsilon
}e^{\left| x\right| (h(x/\left| x\right| )+\gd(\varepsilon) )}\,,
\]
where $\lim_{\varepsilon\to 0}\gd(\varepsilon)=0$, which gives
(\ref{estimate}).

Suppose now that $u$ satisfies (\ref{estimate}). We need to prove
that $u$ can be represented as in (\ref{represent}). Let
$G_{\,\varepsilon }$ be the $\varepsilon $-neighborhood of $G$ and
$h_{\,\varepsilon }=h+\varepsilon$ be the indicator function of
$\overline{G_{\,\varepsilon }}$. Consider the following
Fr\'{e}chet spaces of test functions:
$$ W_{m,\varepsilon }=\left\{ \phi \in H_{loc}^m(\Real^n)\,|\,
<\phi>_{m,\delta}\, <\infty ,
\,\forall\; 0<\delta<\varepsilon \right\},$$ where
$$<\phi>_{m,\delta}\,: =\sup_{\gamma \in \Gamma }\left\{\left|
\left| \phi \right|\right|_{H^m(K+\gamma)}e^{(h_{\delta}(\gamma
/\left| \gamma \right| )\left| \gamma \right| )}\right\} \, .$$
It is obvious that the operator $L^{*}$ maps continuously
$W_{2,\varepsilon } $ into $W_{0,\varepsilon }$. Consequently, the
linear functional
\[
<u,\phi >:=\int_{\mathbb{R}^n} u(x)\phi (x)dx
\]
is continuous on the space $W_{0,\varepsilon }$ for any
$\varepsilon >0$. Since $Lu=0$, Schauder elliptic estimates
together with the periodicity of the operator show that estimates
similar to (\ref{estimate}) hold also for the
derivatives of $u$. One observes that $u$ is a continuous functional on $%
W_{0,\varepsilon }$ which annihilates the range of the operator $%
L^{*}:W_{2,\varepsilon }\rightarrow W_{0,\varepsilon }$. Now
Floquet theory arguments analogous to the ones used in \cite[Section 3.2]{Ku}
can be applied to yield (\ref{represent}). Let us make this part
more precise.

Our first goal is to obtain a Paley-Wiener type theorem for the
Floquet transform in the spaces $W_{m,\varepsilon }$. Let us
denote by $V_{\,\varepsilon }$ the domain in $(\CC^{\, *})^n$
 \bea V_{\,\varepsilon }=\left\{
z=(z_1,\ldots ,z_n)\in (\CC^{\, *})^n\,|\,\right. && \hspace{-5mm}
z_j =\exp ik_j \mbox{ such that }\nonumber\\ &&\hspace{-8mm} Im \,
k=(Im \, k_1,\dots ,Im \, k_n)\in (-G_{\,\varepsilon
})\left.\right\} .\nonumber \eea
 and let
 \bea V^*_{\,\varepsilon}=\left\{ z=(z_1,\ldots ,z_n)\in (\CC^{\, *})^n
 \,|\,z^{-1}=(z_1^{-1},\ldots,z_n^{-1}) \in V_{\,\varepsilon
}\right\} , \nonumber \eea
 The domains $V_{\,\varepsilon }$
form a basis of neighborhoods of the tube $V$, where $V$ is
defined by (\ref{tubulard}). The following statement is a
Paley-Wiener type theorem for the transform ${\cal U}$ which is
suitable for our purpose.
\begin{lemma}\label{P-W}
\begin{enumerate}
\item  The operator
\[
{\cal U}:W_{m,\varepsilon }\rightarrow \Gamma (V^{*}_{\,\varepsilon
},{\cal E}_m)
\]
is an isomorphism, where $\Gamma (V^{*}_{\,\varepsilon },{\cal
E}_m)$ is the
space of holomorphic sections over $V^{*}_{\,\varepsilon }$ of the bundle
${\cal E}_m$, equipped with the topology of uniform convergence on compacta.

\item Under the transform ${\cal U}$, the operator
\[
L^{*}:W_{2,\varepsilon }\rightarrow W_{0,\varepsilon }
\]
becomes the operator ${\cal L}(z)$ of
multiplication by a holomorphic
Fredholm morphism between the fiber bundles ${\cal E}_2$ and
${\cal E}_0$:
\[
\Gamma (V^{*}_{\,\varepsilon },{\cal E}_2)\stackrel{{\cal L}(z)}{\rightarrow }%
\Gamma (V^{*}_{\,\varepsilon },{\cal E}_0).
\]
Here ${\cal L}(z)$ acts on each fiber of ${\cal E}_2$ as the
restriction to this fiber of the operator $ L^{*}$ acting
between $H^2(K)$ and $L^2(K)$.
\end{enumerate}
\end{lemma}

Let us choose a value $\varepsilon _0>0$ such that the
intersection of $\Phi $ with $V_{\,\varepsilon }$ is smooth and
connected. This is possible according to Lemma
\ref{Lambda-lemma2}. From now on, we will only consider the values
$0<\varepsilon <\varepsilon _0$.

Since the image ${\cal U}u$ of the solution $u$ under the Floquet transform
${\cal U}$ is a continuous linear functional on $\Gamma (V^{*}_{\,\varepsilon },%
{\cal E}_0)$ which is in the cokernel of the operator
\[
\Gamma (V^{*}_{\,\varepsilon },{\cal E}_2)\stackrel{{\cal L}(z)}{\rightarrow }%
\Gamma (V^{*}_{\,\varepsilon },{\cal E}_0),
\]
our task is to describe all such functionals. Several
theorems of this kind were proven in \cite{Ku}. In our current
situation such a representation can be obtained rather easily, due
to the simplicity of the structure of the Floquet variety inside
$V_{\,\varepsilon }$. Namely, let $u_z(\cdot)=z^xp(z,\cdot)$ be
the Bloch solution of the equation $Lu=0$ introduced in Lemma
\ref{Lambda-lemma2}. Let also ${\cal H}(\Phi _{\,\varepsilon })$
be the space of holomorphic functions on $\Phi _{\,\varepsilon
}=\Phi \cap V_{\,\varepsilon }$ equipped with the topology of
uniform convergence on compacta. We introduce the mapping
\[
t:\Gamma (V^{*}_{\,\varepsilon },{\cal E}_0)\rightarrow {\cal
H}(\Phi _{\,\varepsilon })
\]
which for any section $f(z,x)$ of the bundle ${\cal E}_0$ produces
\[
t_{f}(z)=<f(z^{-1},\cdot),u_{z}>=\int\limits_{\mathbb{T}^n}f(z^{-1},x)u_{z}(x)dx.
\]
Here $z^{-1}=(z_1^{-1},\ldots, z_n^{-1})$.

\begin{lemma} \label{Isomorphism} Let $0<\varepsilon <\varepsilon
_0$, where $\varepsilon _0$ is the value defined above. Then the
mapping $t$ is a topological homomorphism and the following
sequence is exact: $$
\Gamma (V^{*}_{\,\varepsilon },{\cal E}_2)\stackrel{{\cal L}(z)}{\rightarrow }%
\Gamma (V^{*}_{\,\varepsilon },{\cal E}_0)\stackrel{t}{\rightarrow }{\cal H}%
(\Phi _{\,\varepsilon })\rightarrow 0.
\label{sequ}
$$
\end{lemma}

This lemma practically finishes the proof of the theorem. Namely,
the solution $u$ after the Floquet transform leads to a continuous
linear functional on $\Gamma (V^{*}_{\varepsilon },{\cal E}_0)$
that annihilates the range of the operator of multiplication by
${\cal L}(z)$. Lemma \ref{Isomorphism} implies that such a
functional can be pushed down to the space ${\cal H}(\Phi
_{\,\varepsilon })$. Since this functional, due to the estimate
(\ref{estimate}), is continuously extendable to ${\cal H}(\Phi
_{\,\varepsilon })$ for arbitrarily small values of $\varepsilon
$, it is in fact a hyperfunction (analytic functional) $\mu $ on
$\Phi =\bigcap\limits_{\varepsilon
>0}\Phi _{\,\varepsilon }$. Hence, the action  $<u,\phi >$ of the functional
$u$ on a function $\phi\in W_{0,\varepsilon }$ can be obtained as
\[
<u,\phi >=<\mu (z),t(z)({\cal U}\phi )>.
\]
Applying now the explicit formulas for the transforms ${\cal U}$
and $t$, one arrives to the representation (\ref{represent}). Indeed,
\begin{eqnarray} \label{range}
t_{({\cal U}\phi )}(z)=\int\limits_{K}{\cal U}\phi
(z^{-1},x)u_{z}(x)dx  \\ =\sum \limits_{\gamma \in
\Gamma}\int\limits_{K-\gamma } \phi
(x)z^{-\gamma}u_{z}(x+\gamma)dx \nonumber \\ =\int\limits_{\Real^n
}\phi (x)u_{z}(x)dx. \nonumber
\end{eqnarray}
In this calculation we used the property of the Bloch solutions $$
u_{z}(x+\gamma)=z^{\gamma}u_{z}(x). $$ Therefore, $$
<u,\phi >=<<\mu (z),u_{z}>,\phi >,
$$
which concludes the proof of the theorem. \qed

\mysection{Liouville-type theorem}\label{Sec4}

In this section we discuss Liouville theorems for periodic
equations. We will consider at the moment an arbitrary linear
elliptic operator $P(x,D)$ with smooth $\Gamma $-periodic
coefficients which satisfies the assumptions made in Section
\ref{notations} (as above, without loss of generality we can
reduce the consideration to the case $\Gamma =\ZZ^n$).

\begin{definition} \label{defLiouv} {\em We say that {\em the Liouville theorem
holds true for the operator $P$}, if for any $N\in \mathbb{N}$ the
space $V_{N}(P)$ of solutions of the equation $Pu=0$ in $\Real^n$
that can be estimated as $$ \left| \left| u\right| \right|
_{L^2(K+\gamma )} \leq C(1+\left| \gamma \right| )^N \mbox{ for
all } \gamma \in \Gamma $$ is finite dimensional.}
\end{definition}

In the case when the Liouville theorem holds, we will be also
interested in the dimensions  $d_N$ of the spaces $V_{N}(P)$ and
in representations of their elements analogous to (\ref{polyn}).

The result below explains under what
conditions on the operator $P$ a Liouville-type theorem holds. These
conditions will then be verified for some specific classes of
operators.

As was mentioned in the introduction, solutions representable as
(\ref{polyn}) are just Floquet solutions with zero quasimomentum.
So, the Liouville theorem of \cite{AL,MS} cited in the
introduction states that any polynomially growing solution is a
Floquet solution with a zero quasimomentum. Let us also mention
that any Bloch solution $e^{ik \cdot x}p(x)$ with a real
quasimomentum $k$ is automatically bounded. This means that the
validity of the Liouville theorem for an operator $P$ implies that
the number of the real quasimomenta of solutions of the equation
$Pu=0$ must be finite (modulo the action of the reciprocal
lattice). In other words, the Fermi surface for $P$ intersects the
real space at a finite number of points (modulo the reciprocal
lattice). In terms of the Floquet variety it means that the set
$\mathcal{Z}:=\Phi_P \cap T$ is finite. We denote the cardinality
of a set $A$ by $\#A$. As the second statement of the next theorem
shows, the finiteness of $\mathcal{Z}$ is in fact the only claim
of the Liouville theorem.

\begin{theorem}
\label{Liouville}
\begin{enumerate}
\item The equation $Pu=0$ has a nonzero polynomially growing
solution if and only if it has a nonzero bounded Bloch solution,
i.e. if and only if the intersection $F_{L} \cap \Real^{n}$ of the
Fermi surface for $P$ with the real space is not empty (or
equivalently, $\mathcal{Z}=\Phi_P \cap T \neq \emptyset$).

\item The Liouville theorem holds for the operator $P$
if and only if the intersection $F_P \cap \Real^n$  is a finite
set modulo the reciprocal lattice (or equivalently,
$\#\mathcal{Z}<\infty$). Moreover, if $\#\mathcal{Z}=\infty$  then
the Liouville theorem does not hold even for bounded solutions,
i.e., $d_0=dim(V_{0})=\infty$.

\item  If the Liouville theorem holds,
then each solution $u \in V_N(P)$ can be represented as a finite
sum of Floquet solutions:

\begin{equation}
u(x)= \sum\limits_{q \in F_P \cap \Real^n} e^{iq \cdot
x}\sum\limits_{|j|\leq N}x^{j} p_{j,q}(x). \label{li_repres}
\end{equation}

\item If the Liouville theorem holds,
then for all $N\geq 0$ we have $$ d_N\leq d_0 q_{n,N}<\infty\,, $$
where $q_{n,N}$ is the dimension of the space of all polynomials
of degree at most $N$ in $n$ variables.

\item Assume that the Liouville theorem holds and that for each
real quasimomentum $q$ (i.e., for each $q\in F_P \cap \Real^n$)
the conditions of Theorem \ref{Fl_dimen} are satisfied. Then for
each $N \geq 0$ the dimension $d_N$ of the space $V_N(P)$ is equal
to the sum over $q \in (F_P \cap \Real^n)/\Gg^*$ of the dimensions
of the spaces of $\lambda_q$-harmonic polynomials (see Definition
\ref{Qharm}), where $\lambda_q$ is the first nonzero homogeneous
term in the Taylor expansion at the point $q$ of the dispersion
relation (band function) $\lambda(k)$.
\end{enumerate}
\end{theorem}

\pf Statements $4$ and $5$ follow from $3$ together with Lemma
\ref{Floq_struct} and Theorem \ref{Fl_dimen}. So, we first prove
statements $2$ and $3$ and conclude with the proof of the first
statement.

In order to prove $2$ let us notice that if $\#\mathcal{Z}=\infty$
then each point $z=\exp ik\in \mathcal{Z}$ provides a bounded
Bloch solution with the quasimomentum $k$, and these solutions are
linearly independent. This means that the Liouville theorem cannot
hold in this case.

Assume now that $\# \mathcal{Z} < \infty$. We need to prove that
the Liouville theorem and representation (\ref{li_repres}) hold
true. Obviously, if $u$ has a representation of the form
(\ref{li_repres}), then $u$ is of a polynomial growth. The proof
that any polynomially growing solution is of the form
(\ref{li_repres}) follows the same simple strategy as in the
proofs of Theorem \ref{hyperfunction} and as is in the proof of
the main Floquet representation \cite[Theorem 3.2.1]{Ku} (which,
in turn, comes from the approach of \cite{E} and \cite{Pa1}). As in
the case with the ``fundamental principle'' (see \cite{E} and
\cite{Pa1}), it is more convenient to deal with a dual
formulation, as it is done in \cite {Ku}. Namely, any polynomially
growing solution $u(x)$ can be interpreted in the dual way, as a
functional on an appropriate functional space, which belongs to
the cokernel of the dual operator $P^{*}$. Consequently, a
representation theorem for all such functionals must be obtained.
In order to make this idea precise, we need to introduce
appropriate test functions spaces.

Consider the Fr\'{e}chet spaces
\[
C_m=\left\{ \phi \in H_{loc}^m(\Real^n)|\;\sup_{\gamma \in \Gamma
}\left| \left| \phi \right| \right| _{H^m(K+\gamma )}(1+\left|
\gamma \right|)^N<\infty ,\;\forall N\right\} .
\]
Let the order of the operator $P$ be $m$, then it is clear that $P^{*}$ maps
continuously $C_m$ into $C_0$. Due to the polynomial growth of $u(x)$, the
linear functional
\[
<u,\phi >=\int_{\mathbb{R}^n} u(x)\phi (x)dx
\]
is continuous on $C_0$. Since $Pu=0$, one easily observes that $u$
annihilates the range of the operator $P^{*}:C_m\rightarrow C_0$. We need
now a Paley-Wiener type theorem for the spaces $C_m$ with respect to the
Floquet transform.

\begin{lemma}
\label{P-WC}
\begin{enumerate}
\item  The operator
\[
{\cal U}:C_m\rightarrow C^\infty (T,{\cal E}_m)
\]
is an isomorphism, where $C^\infty (T,{\cal E}_m)$ is the space of
$C^\infty $ sections of the bundle ${\cal E}_m$ over the complex
torus $T$, equipped with the standard topology.

\item  Under the transform ${\cal U}$, the operator
\[
P^{*}:C_m\rightarrow C_0
\]
becomes the operator ${\cal P}(z)$ of
multiplication by a holomorphic Fredholm morphism between the
fiber bundles ${\cal E}_m$ and ${\cal E}_0$:
\[
C^\infty (T,{\cal E}_m)\stackrel{{\cal P}(z)}{\rightarrow }C^\infty (T,{\cal %
E}_0).
\]
Here ${\cal P}(z)$ acts on each fiber of ${\cal E}_m$ as the
restriction  to this fiber of the operator $ P^{*}$ acting between
$H^m(K)$ and $L^2(K)$.

\item  The operator ${\cal P}(z)$ is invertible for a point $z\in T$ if and
only if $z^{-1} \notin \Phi $.
\end{enumerate}
\end{lemma}

The next lemma is an analog of the classical theorem on the
structure of distributions supported at a single point. Together
with the previous lemma it essentially leads to the statement of
the theorem.
\begin{lemma}
\label{delta} Let $T$ be a $C^\infty $-manifold and ${\cal
P}:T\rightarrow L(B_1,B_2)$ be a $C^\infty $-function with values
in the space $L(B_1,B_2)$ of bounded linear operators between
Banach spaces $B_1$ and
$B_2$. Assume that for each  $z\in T$ the operator ${\cal P}(z)$
is a Fredholm operator. Then

\begin{enumerate}
\item  If ${\cal P}(z)$ is surjective for all points $z$ in $T$, then the
multiplication operator
\[
C^\infty (T,B_1)\stackrel{{\cal P}(z)}{\rightarrow }C^\infty (T,B_2)
\]
is surjective.

\item  If ${\cal P}(z)$ is surjective for all points $z$ except
a finite subset ${\cal Z} \subset T$,
then any continuous linear functional $g$ on the space of smooth vector
functions $C^\infty (T,B_2)$ that annihilates the range of the
multiplication operator
\[
C^\infty (T,B_1)\stackrel{{\cal P}(z)}{\rightarrow }C^\infty (T,B_2)
\]
has the form
\begin{equation}
<g,\phi >=\sum\limits_{z \in {\cal Z}}\left[ \sum\limits_{j\leq N}
D_{j,z}(<g_{j,z},\phi >)\right] _{z}.  \label{repr}
\end{equation}
Here $g_{j.z}$ are continuous linear functionals on $B_2$,
$<g_{j,z},\phi>$ denotes the duality between $B_2^{*}$ and $B_2$,
$D_{j,z}$ are
linear differential operators on $T$, and $N\in \mathbb{N}$.
\end{enumerate}
\end{lemma}

We are ready now to finish the proof of the nontrivial part of the
third statement of Theorem \ref{Liouville}.

If $u$ is a solution of polynomial growth, it belongs, as it has
been mentioned already, to the cokernel of the operator
$P^{*}:C_m\rightarrow C_0$. After the Floquet transform we are
dealing with the cokernel of the operator
\[
C^\infty (T,{\cal E}_2)\stackrel{{\cal P}(z)}{\rightarrow }C^\infty (T,{\cal %
E}_0).
\]
By Lemma \ref{P-WC}, the only points $z\in T$ where ${\cal P}(z)$
is not invertible are those points where $z^{-1}$ belongs to the
Floquet variety. Since by our assumption the set ${\cal Z}=T\cap
\Phi$ is finite, it follows that the operator function ${\cal
P}(z)$ satisfies all the assumptions of Lemma \ref{delta}. The
fact that we are dealing with Banach bundles instead of fixed
Banach spaces is irrelevant, since these bundles are trivial. This
means that we have the representation (\ref{repr}) with $g_j\in
L^2(\mathbb{T}^n )$. According to Lemma \ref{Floq_struct},
functionals of the form (\ref{repr}) correspond under the inverse
Floquet transform exactly to functions of the form (\ref{polyn}).

It remains to prove the first statement of the theorem.  Let $u$
be a polynomially growing solution. Assume that
$\mathcal{Z}=\emptyset$, i.e., the intersection of the Floquet
variety $\Phi $ with the complex torus $T$ is empty. Therefore,
the last statement of Lemma \ref{P-WC} implies the invertibility
of $P(z)$ for all $z\in T$. Now, the first statement of Lemma
\ref{delta} guarantees the surjectivity of the mapping
\[
C^\infty (T,{\cal E}_m)\stackrel{{\cal P}(z)}{\rightarrow }C^\infty (T,{\cal %
E}_0)
\]
and hence the absence of any nontrivial functionals on $C^\infty
(T,{\cal E}_0)$ that annihilate the image of this mapping. Since
under the Floquet transform ${\cal U}$, a polynomially growing
solution $u(x)$ is mapped to such a functional, we conclude that
$u=0$. \qed
\begin{remark}
\label{alterproof} {\em The first statement of Theorem
\ref{Liouville} is a part of the analog of the Bloch theorem
provided in Theorem 4.3.1 of \cite{Ku}. Namely, the existence of a
sub-exponentially (in particular, polynomially) growing solution
implies the existence of a Bloch solution with a real
quasimomentum, and hence the nonemptiness of the real Fermi
variety. For completeness, we gave above an independent proof of
this statement.}
\end{remark}

One realizes now that the cases when a Liouville-type theorem
holds in a nonvacuous way are extremely rare. Namely, Theorem
\ref{Liouville} shows that this happens only when the Fermi
variety touches the real subspace at a finite set of points
(modulo the reciprocal lattice). This means in particular, that in
the {\em selfadjoint case}, one should expect this to happen only
at the edges of the spectral gaps. Although it is possible to
imagine interior points of the spectrum where such a thing could
occur, it is hard to believe that these cases could be anything
more than accidents.

One can expect the following conjecture to be true:
\begin{conjecture}
Let $P$ be a ``generic" self-adjoint second order elliptic
operator with periodic coefficients and
$(\lambda_{-},\lambda_{+})$ be a nontrivial gap in its spectrum.
Then each of the gap's endpoints is a unique (modulo the dual
lattice) and nondegenerate extremum of a single band function
$\lambda_{j}(k)$.
\end{conjecture}

The validity of this conjecture together with Theorem
\ref{Liouville} would imply that generically at the gap ends the
dimension of the space $V_{N}$ is equal to the dimension $h_{n,N}$
of the space of all harmonic polynomials of order at most $N$ in
$n$ variables. Unfortunately, the only known theorem of this kind
is the recent result of \cite{KR}, which states that generically a
gap edge is an extremum of a single band function.

At the bottom of the spectrum, however, much more is known. The
theorem below combines some results of \cite{FKT, KS, Pr} with the
statement of Theorem \ref{Liouville} to obtain the structure and
dimension of the space of polynomially growing solutions in this
case. Below the spectrum, the Liouville theorem holds vacuously,
according to the first statement of Theorem \ref{Liouville} and
Theorem 5.5.1 in \cite{Ku}.

\begin{theorem}
\label{applicat}
\begin{enumerate}

\item Let $H=-\Delta + V(x)$ be a Schr\"{o}dinger operator with a periodic
real valued potential $V \in L^{r/2}(\mathbb{T}^{n}),\,r>n$. Then
the lowest band function $\lambda_1(k)$ has a unique nondegenerate
minimum $\Lambda_{0}$ at $k=0$. All other band functions are
strictly greater than $\Lambda_{0}$. Every solution $u\in
V_N(H-\Lambda_{0})$ is representable in the form (\ref{polyn}).
The dimension of the space $V_N(H-\Lambda_{0})$ is equal to
$h_{n,N}$.

\item Let $V$ be like in the previous statement, then there exists $\epsilon>0$
such that for any periodic real valued magnetic potential $A$ such
that $$ \left|\left| A \right|\right|_{L^{r}(\mathbb{T}^{n})} <
\epsilon $$ and
\begin{equation}\label{A0}
\int\limits_{\mathbb{T}^n}A(x)dx=0 \end{equation}
 the following statements hold true: The
lowest band function $\lambda_1(k)$ of the magnetic
Schr\"{o}dinger operator $H=(i \nabla +A)^2+V$ attains a unique
nondegenerate minimum $\Lambda_{0}$ at a point $k_{0}$. All other
band functions are strictly greater than $\Lambda_{0}$. Every
solution $u\in V_N(H-\Lambda_{0})$ is representable in the Floquet
form $$ v(x)=e^{ik_{0}\cdot x}\sum\limits_{|j|\leq N}x^jp_j(x) $$
with periodic functions $p_j(x)$. The dimension of the space
$V_N(H-\Lambda_{0})$ is equal to $h_{n,N}$.

\item Suppose that $L$ is a second order elliptic operator of the form
(\ref{operator}) such that $\Lambda_{0} \geq 0$.

If $\Gl(0)=0$ (i.e. $0\in \Xi$), then the Liouville theorem holds
and every solution $u\in V_N(L)$ is representable in the form
(\ref{polyn}). The dimension of the space $V_N(L)$ is equal to
$h_{n,N}$ in the case when $\Lambda_{0}=0$, and to $q_{n-1,N}$
when $\Lambda_{0}>0$.

If $\Gl(0)>0$ then the equation $Lu=0$ does not admit a nontrivial
polynomially growing solution. So, the Liouville theorem holds
vacuously.
\end{enumerate}
\end{theorem}
\pf  1. The result of \cite{KS} says that the lowest band function
$\lambda_1(k)$ has a unique nondegenerate minimum $\Lambda_{0}$ at
$k=0$ and that all other band functions are strictly greater than
$\Lambda_{0}$. Now Theorem \ref{Liouville} implies the rest of the
claims of this statement. \vspace{2mm}

\noindent 2. When both the electric and magnetic potentials are
sufficiently small, then the result of \cite{FKT} states that the
lowest band function $\lambda_1(k)$ of the magnetic
Schr\"{o}dinger operator $H=(i \nabla +A)^2+V$ attains a unique
nondegenerate minimum $\Lambda_{0}$ at a point $k_{0}$, while all
other band functions are strictly greater than $\Lambda_{0}$. This
statement, however, can be easily extended to the case of
arbitrary electric and small magnetic potential. Indeed, when the
magnetic potential is equal to zero, one can refer, as in the
previous case, to \cite{KS}. At this moment one has to use
analyticity of the Bloch variety. Namely, the statement of Lemma
\ref{Fermi} (see also \cite[Theorem 4.4.2]{Ku}) can be easily
extended to include analyticity with respect to the potentials
(see, for instance, \cite{FKT}). More precisely, there exists an
entire function $f(k,\lambda,A,V)$ of all its arguments such that
$f(k,\lambda,A,V)=0$ is equivalent to
 $$ (k,\lambda)\in B_{(i\nabla +A)^2+V}\,, $$
 where $B_H$ is the Bloch variety of the operator $H$.
 Now, the result of \cite{KS} for $A=0$
together with the stated analyticity property imply the required
features of the lowest band function for sufficiently small
magnetic potentials. The last step is to use  again Theorem
\ref{Liouville}. Note that the normalization (\ref{A0}) always can
be achieved by a gauge transformation which does not affect the
spectrum and the Liouville property. \vspace{2mm}

\noindent 3. The assumption $\Lambda(0)\geq 0$ implies that the
operator $L$ admits a positive periodic supersolution. It follows
from Lemma \ref{constants} that the Fermi surface $F_L$ can touch
the real space only at the origin (modulo the reciprocal lattice
$\Gamma ^{*}$) and in this case $\Lambda(0)= 0$.  Therefore, by
the first part of Theorem \ref{Liouville}, the Liouville Theorem
holds vacuously if $\Lambda(0)> 0$.

Suppose now that $\Lambda(0)= 0$. Lemma \ref{Lambda-lemma} implies
that if $\Lambda_{0}>0$ then the point $k=0$ is a noncritical
point of the dispersion relation, and if $\Lambda_{0}=0$ then
$k=0$ is a nondegenerate extremum. Now Theorem \ref{Liouville}, as
before, completes the proof. \qed

\mysection{Proofs of the lemmas}\label{Sec5}

\noindent \mbox{{\bf Proof of Lemma \ref{Floq_struct}}: } The
first claim of the lemma corresponds to Theorem 3.1.3 in
\cite{Ku}. In order to prove the second part of the lemma, let us
fix a $k_{0}\in \mathbb{C}^n$, and choose a closed subspace
$M\subset H^{m}(\mathbb{T}^n)$ complementary to the kernel of the
operator $P^*(x,D-k_{0})$. Consider the (analytically depending on
$k$ in a neighborhood of $k_{0}$) subspace $$
\Pi(k):=P^*(x,D-k)(M)\subset L^2(\mathbb{T}^n). $$ and $$
\mathcal{N}:=\left[ \Pi(k_{0}) \right]^{\perp}. $$ Then
$dim(\mathcal{N})=a_{k_{0}}$, and for values of $k$ close to
$k_{0}$ the space $\mathcal{N}$ remains a complementary subspace
to $\Pi(k)$. Representing the operators $P^*(x,D-k)$ in the matrix
form according to the decompositions $$
H^{m}(\mathbb{T}^n)=M\oplus Ker\,P^*(x,D-k_0) $$ and $$
L^2(\mathbb{T}^n)=\Pi(k)\oplus \mathcal{N}, $$ we get the matrix
$$ \left(
\begin{array}{cc}
B(k)&* \\ 0&C(k)
\end{array}
\right), $$ where $B(k)$ is an invertible analytic operator
function, and $C(k)$ is an analytic matrix function of the size
$a_{k_0}\times a_{k_0}^{*}$. Here $a_{k_0}^{*}$ is the dimension
of the kernel of the operator $P^*(x,D-k_{0})$. (Notice that
$a_{k_0}=a_{k_0}^{*}$ if $ind\,P=0$, which is true for instance,
when dealing with scalar elliptic operators, due to the
Atiyah-Singer theorem.) Now, the space of all distributions
orthogonal to the range of $P^*$ and supported at $\exp(-ik_{0})$
reduces to the space of all distributions supported at $k_{0}$,
acting on $\CC^{\,a_k}$-valued vector functions, and orthogonal to
the range of the operator of multiplication by $C(k)$. If we drop
the orthogonality condition, the dimension of the space of all
such distributions of order at most $N$ is obviously equal to
$a_{k}q_{n,N}$, which proves the estimate. We point out that a
direct proof of this estimate for scalar operators can be also
easily derived using the Leibnitz's rule. \qed

\noindent \mbox{{\bf Proof of Lemma \ref{Lambda-lemma}}: }
Statements 1 through 3 of the lemma are contained in \cite{LP},
except the statement that the geometric rather than the algebraic
multiplicity of the eigenvalue $\Lambda(\xi)$ is equal to one. The
latter follows easily from Lemma 5.2 of \cite{LP}. Alternatively,
it can be deduced from general theorems on positive operators
defined on an ordered Banach space (see for instance,
\cite[Theorem 2.10]{Kr}). Statement 4 is proven in \cite[Theorem
5]{Pr}. \qed

\noindent \mbox{{\bf Proof of Lemma \ref{Lambdazero}}: }
Statements 1--3 follow from the results of \cite{A1a, LP}, while
statements 4--5 follow from \cite[Theorem 5]{Pr}. \qed

\noindent \mbox{{\bf Proof of Lemma \ref{analytic}}: } Consider
the following family of operators on the torus:
$L(x,D-i\xi)-\Lambda(\xi)$. It follows from Lemma
\ref{Lambda-lemma} that this family is analytic in a complex
neighborhood $W$ of the set $G$ and its values are Fredholm
operators between the appropriate Sobolev spaces. The same lemma
implies that the dimension of the kernel of all these operators is
equal to $1$. Hence, these kernels form an analytic fiber bundle
over $W$ (see Theorem 1.6.13 and the corresponding references in
\cite{Ku}). One can always assume that the domain $W$ is convex
(in the geometric sense). Then the kernel bundle (as all vector
bundles on $W$) is topologically trivial. Since $W$, being convex,
is a domain of holomorphy (see for instance Corollary 2.5.6 in
\cite{Ho}), therefore, the result of \cite{Gr} (an instance of the
so called Oka's principle) implies that the bundle is also
analytically trivial. This means the existence of a nowhere zero
analytic section $u_\xi$. Positivity of $u_\xi$ for $\xi \in \Xi$
can be achieved as follows. Let us choose any nonzero analytic
solution $u_\xi$ as above. Then for some small neighborhood
$W_1\subset W$ of $G$, we have $u_\xi(0)\neq 0$. So, we may
normalize $u_\xi$ by dividing it by $u_\xi(0)$. The resulting
solution is clearly positive for $\xi \in \Xi$. \qed

\noindent \mbox{{\bf Proof of Lemma \ref{Fermi-lemma}}: } 1. Let
$u(x)=e^{ik\cdot x}p(x)$ be a nonzero Bloch solution, where $p(x)$
is a $\Gamma $-periodic function, and $k\in F\cap {\cal T}$.
Assume first that $Im \, k\in -\stackrel{\circ}{G}$, so, $\Gl_0
>0$. We need to prove that $Re \, u=Im \, u=0$. We show for
instance, that $u_1:=Re \, u=0$. Suppose that $u_1\neq 0$. We may
assume that $u_1(x_1)>0$, for some $x_1\in \Real^n$. Consider the
positive solution
 $$v(x)=\int\limits_\Xi u_{\,\xi}(x)\,d\gs(\xi),$$
  where $d\gs$ is the $(n-1)$-dimensional surface area element on $\Xi$.
For every $M>0$ there exists $R>0$ such that $v(x)-Mu_1(x) >0$ for
all $|x|>R$. By the generalized maximum principle, $v(x)>Mu_1(x)$
in $\mathbb{R}^n$. Since $M$ is arbitrarily large and
$u_1(x_1)>0$, we arrived at a contradiction. Note that this
argument applies also to {\em any Floquet} solution with a
quasimomentum $k$ such that $Im \, k\in -\stackrel{\circ}{G}$.

Suppose now that $Im \, k\in - \Xi$ and $\Gl\geq 0$. Clearly, it
is enough to show that there exists a real constant $C$ and $\xi
\in \Xi$ such that $u_1:=Re \, u=C u_{\,\xi}$. Let $\xi=-Im \, k$.
Then for a sufficiently small $\varepsilon>0$ the function
$v_\varepsilon :=\frac{u_{\,\xi}}{2}-\varepsilon u_1$ is a
positive solution of the equation $Lu=0$, which is smaller than
$u_{\,\xi}$. Recall that $u_{\,\xi}$ is a minimal positive
solution of the equation $Lu=0$. Therefore, there exists $c>0$
such that $v_\varepsilon=c u_{\,\xi}$, which implies that $u_1=C
u_{\,\xi}$ for some $C\in \Real$. \vspace{2mm}

\noindent 2.  Consider the zero set $F_1$ of the analytic function
$\Lambda(ik)$ in a small complex neighborhood of $-i\Xi$. Since
$\Lambda_0 >0$, it follows that the gradient of $\Lambda(ik)$ is
not zero on $-i\Xi$. Therefore, $F_1$ is a smooth analytic
variety. We will show that the Fermi surface $F$ coincides with
$F_1$ in a neighborhood of $-i\Xi$, which will conclude the proof
of the lemma. Indeed, obviously $F_1 \subset F$. Consider a point
$k_0=-i\xi_0\in -i\Xi$. By Lemma \ref{Lambda-lemma}, zero is a
simple eigenvalue of the operator $L(x,D+k_0)=L(x,D-i\xi_0)$. This
means that the spectral projector that corresponds to a
neighborhood of zero is one-dimensional for all complex $k$ close
to $k_0$. We conclude that for all $k$ in a complex neighborhood
of $-i\Xi$ there is exactly one eigenvalue close to zero of the
operator $L(x,D+k)$. By Lemma \ref{analytic}, we know this
eigenvalue, namely $\Lambda(ik)$. Let now $k$ belongs to a small
neighborhood of $-i\Xi$ and assume that $k\not\in F_1$. Then
$\Lambda(ik) \neq 0$, and hence zero cannot be the eigenvalue of
$L(x,D+k)$. This means that $k$ does not belong to the Fermi
surface $F$. \qed

\noindent \mbox{{\bf Proof of Lemma \ref{constants}}: }
 1. If $c\gneqq 0$, the assertion of the lemma follows from
\cite[Theorem 4.5]{LP}. On the other hand, if $c=0$, then $0\in
-i\Xi$, and in particular,  $0\in (-G)$. It follows from Lemma
\ref{Fermi-lemma} that any Bloch solution with a real
quasimomentum is the constant solution. \vspace{2mm}

\noindent 2. This assertion follows directly from the part 1 using
the operator $\psi^{-1}L\psi$. \qed

\noindent \mbox{{\bf Proof of Lemma \ref{P-W}}: } The second
statement of the lemma coincides with Theorem 2.2.3 in \cite{Ku}.
So, we need to prove only the first statement.

Let $\varphi \in W_{m,\varepsilon }$. We will show that the series
(\ref{Floquet}) converges uniformly on compacta in
$V^*_{\,\varepsilon }$ as a series of functions on
$V^*_{\,\varepsilon }$ with values in $H^m(K)$. This would imply
that $\mathcal{U}\varphi \in \Gamma (V^*_{\,\varepsilon
},H^m(K))$, and that the corresponding (one-to-one) mapping $\mathcal{U}%
:W_{m,\varepsilon }\rightarrow \Gamma (V^*_{\,\varepsilon },H^m(K))$ is
continuous. Let $0<\delta <\delta _1<\varepsilon $. Let $z=\exp
ik \in V^*_{\,\delta }$ which means that, $Im\, k\in
G_{\,\delta}$. We have

\begin{eqnarray}
\left| \left| \mathcal{U}\varphi (z,\cdot )\right| \right|
_{H^m(K)}\leq \sum\limits_{\gamma \in \Gamma }\left| \left|
\varphi \right| \right| _{H^m(K-\gamma )}e^{-Im\,k\cdot \gamma }=
\sum\limits_{\gamma \in \Gamma }\left| \left| \varphi \right|
\right| _{H^m(K+\gamma )}e^{Im\,k\cdot \gamma } \nonumber \\ \leq
\sum\limits_{\gamma \in \Gamma }\left| \left| \varphi \right|
\right| _{H^m(K+\gamma )}e^{(h(\gamma /\left| \gamma \right|
)+\delta )\left| \gamma \right| }\leq C_\delta <\varphi
>_{m,\delta_{1}} <\infty . \nonumber
\end{eqnarray}
We need to check now that the mapping $\mathcal{U}$ acts from
$W_{m,\varepsilon }$ into $\Gamma (V^{*}_{\,\varepsilon
},\mathcal{E}_m)$. This amounts to showing that
$\mathcal{U}\varphi $ satisfies the appropriate Floquet boundary
conditions and hence is in fact a section of the sub-bundle
$\mathcal{E}_m\subset V^*_{\,\varepsilon }\times H^m(K)$. This is
a straightforward calculation (see also Theorem 2.2.2 in
\cite{Ku}).

On the other hand, let us assume that $s(z)\in \Gamma (V^*_{\,\varepsilon },%
\mathcal{E}_m)$. If $z=\exp ik$, then $s$ as a function of $k$ is periodic
with respect to the reciprocal lattice $\Gamma ^{*}$. Expanding it into the
Fourier series, we get
\[
s(z)=\sum\limits_{\gamma \in \Gamma }s_{\,\gamma }z^\gamma \text{,}
\]
where $s_{\,\gamma }\in H^m(K)$. We can now define a function $\varphi $ on $%
\Real
^n$ such that $\varphi (x-\gamma )=s_{\,\gamma }(x)$ for $x\in K$ and $%
\gamma \in \Gamma $.

The function $\varphi $ belongs to $H^m$ in the interior of each
of the cubes $K+\gamma $. One only needs to check that it belongs
to $H_{loc}^m$ at the boundary points of these cubes. The
requirement that $s(z)$ is a section of the bundle $\mathcal{E}_m$
rather than just of the bundle $ V^*_{\,\varepsilon }\times
H^m(K)$ does exactly this (see the discussion at the top of page
96 in \cite{Ku}).

It remains to show that $\varphi \in W_{m,\varepsilon }$. We use
the standard formulas for the Fourier coefficients to get
\[
\varphi (\cdot -\gamma )=s_{\,\gamma }=\frac 1{(2\pi
)^n}\int\limits_Bs(e^{i(\gb+i\alpha )})e^{-i(\gb+i\alpha )\cdot
\gamma }\,d\gb,\quad\forall\, \alpha\in G_{\,\varepsilon }\,,
\]
where $B$ is the first Brillouin zone, and we write $z=\exp
ik=\exp (i(\gb +i\ga))$. Note that
\begin{equation}\label{sgest}
\left| \left|\phi\right| \right| _{H^m(K+\gg)} \leq \max_{z\in
V^*_{\,\delta _1}}\left| \left| s(z)\right| \right|
_{H^m(K)}e^{\ga\cdot(-\gamma)}\quad\forall\, \alpha\in
G_{\,\varepsilon }\,,
\end{equation}
and therefore,
\begin{equation}\label{sgest1}
\left| \left|\phi\right| \right| _{H^m(K+\gg)} \leq \max_{z\in
V^*_{\,\delta _1}}\left| \left| s(z)\right| \right|
_{H^m(K)}e^{-(h(\gamma /\left| \gamma \right| )+\delta_{1} )}\,.
\end{equation}
This implies immediately that
\begin{eqnarray}
<\varphi >_{m,\delta }=\sup_{\gamma \in \Gamma }\left\{\left| \left| \varphi
\right| \right| _{H^m(K+\gamma )}e^{(h(\gamma /\left| \gamma
\right| )+\delta )\left| \gamma \right| }\right\}
 \\ \leq
C\max_{z\in V^*_{\,\delta _1}}\left| \left| s(z)\right| \right|
_{H^m(K)}\sup_{\gamma \in \Gamma }e^{-(\delta _1-\delta )\left|
\gamma \right| }<\infty ,\nonumber
\end{eqnarray}
if $\delta _1 > \delta$.
\qed

\noindent \mbox{{\bf Proof of Lemma \ref{Isomorphism}}: } The
statement of this lemma is established in a much more general
situation at the beginning of the proof of Theorem 1.7.1 in
\cite{Ku}. However, for the sake of completeness we provide here
the proof for our simpler particular situation. First of all, the
sequence of the lemma is a complex (i.e., the composition of any
two consecutive operators in it is equal to zero). One needs to
prove this only in the second term of the sequence, where it
follows immediately from the equality (\ref{range}). Indeed, since
$u_z$ solves the equation $Lu=0$, (\ref{range}) followed by
integration by parts proves the statement.

Let us turn to the exactness. We need to prove it in the second
and third terms of the sequence. Consider the second term. Let
$f(z,x)\in \Gamma (V_\varepsilon ^{*},{\cal E}_0)$ be such that
$t_f(z)=0$. This means that for any $z\in \Phi _\varepsilon $ the
function $f(z^{-1},\cdot )$ is orthogonal to the Bloch solution
$u_z$ of the equation $Lu=0$. We need to show that $g(z)={\cal
L}(z^{-1})^{-1}f(z)$ is analytic, which will mean that $f$ belongs
to the range of ${\cal L}$. The function $g(z^{-1})$ is
automatically analytic outside of $\Phi_\epsilon$, so we only need
to make sure that it does not develop any singularities at this
subset. We will show that all the necessary and sufficient
conditions for the analyticity of $g$ have the form of
orthogonality of values of $f$ at certain points to certain
functionals. This would resolve the issue, since all such possible
orthogonality conditions are the orthogonality of $f(z^{-1})$ to
the kernel of $L$ on Bloch functions with a quasimomentum $z$, and
hence to the vanishing of $t_f(z)$. As it was shown in the proof
of Theorem 3.3.1 in \cite[pages 113-114]{Ku}, the inverse operator
to ${\cal L}(z^{-1})$ is the ratio of two analytic functions: $$
{\cal L}(z^{-1})^{-1}=B(z) / \Delta(z),$$ where $B(z)$ is an
analytic function with values in bounded operators from
$L_2({\mathbb T}^n)$ to $H^2({\mathbb T}^n)$, and $\Delta(z)$ is a
scalar analytic function, which is a regularized determinant of
${\cal L}(z^{-1}){\cal L}(z_0^{-1})^{-1}$ for some point $z_0$
where the operator is invertible. Such regularized determinants
are determined in the standard way by the eigenvalues of the
corresponding operators (see for instance Section 2 of Chapter
{\small{\sf IV}} in \cite{GoK} for general definitions and
properties of regularized determinants, and for our particular
situation the proof of Theorem 3.1.7 and related discussion in
Section 1.2 in \cite{Ku}). The simplicity of the eigenvalue
$\Lambda(\xi)$ (Lemma \ref{Lambda-lemma}) implies that if we
introduce instead of $z$ the coordinate $\xi$ such that $z=\exp
\xi$, then $\Delta(\exp\xi)=\Lambda(\xi) \Delta_1(\xi)$, where
$\Delta_1(\xi)$ is an analytic function with no zeros in the
domain under our consideration. We recall now that $\Lambda$ has
simple zeros. Hence, the necessary and sufficient condition for
$f$ to belong to the range of the operator ${\cal L}$ on the space
of analytic sections is that the vector-function $B(z)f(z)$
vanishes on the set of the zeros of $\Lambda$. These conditions
obviously have the form of the orthogonality of values of $f$ to
some functionals. As it was explained above, this implies
exactness at the second term of the sequence.

Let us turn now to proving the exactness at the third term. We
need to show that arbitrary analytic function on $\Phi_\epsilon$
can be obtained as $t_f(z)$ for some $f \in
\Gamma(V^{*}_{\epsilon},\cal{E}_0)$.

Let us denote by $\Phi^{*}_\epsilon$ the manifold $$
\Phi^{*}_\epsilon = \{ z | \,z^{-1} \in \Phi_\epsilon \}. $$
Consider the restriction mapping
\begin{equation}\label{restrict}
 \Gamma(V^{*}_\epsilon, {\cal
E}_0) \rightarrow \Gamma(\Phi^{*}_\epsilon, {\cal E}_0).
\end{equation}
Notice that $\Phi^{*}_\epsilon$ is an analytic subset in
$V^{*}_\epsilon$ and that $V_\epsilon$ and $V^{*}_\epsilon$ are
domains of holomorphy. The latter can be easily proven using power
test functions $z^a$ with integer (but not necessarily
nonnegative) powers $a$ (a similar derivation can be found in the
proof of the implication $(iii) \rightarrow (i)$ of Corollary
2.5.8 in \cite{Ho}). Then Corollary 1 of the Bishop's theorem
\cite[Theorem 3.3]{ZK} (see the original theorem in \cite{Bi})
claims that the restriction mapping (\ref{restrict}) is surjective
(recall that the bundle ${\cal E}_0$ is trivial). Hence, it is
sufficient to prove that the mapping $$ \tilde{t}:
\Gamma(\Phi_\epsilon, {\cal E}_0) \rightarrow {\cal
H}(\Phi_\epsilon). $$ defined as $$
\tilde{t}_{f}(z)=<f(z,\cdot),u_{z}>=\int\limits_{\mathbb{T}^n}f(z,x)u_{z}(x)dx
$$ is surjective. Consider the continuous operator $T(z): L^2(K)
\rightarrow \CC$ defined as
$T(z)y=<y,u_{z}>=\int\limits_{\mathbb{T}^n }y(x)u_{z}(x)dx$. Since
$u_z$ is not zero, this operator is surjective. It is clear that
it depends analytically on $z$. According to Allan's theorem (see
\cite{Al} or Theorem 4.4 in \cite{ZK}), since $\Phi_\epsilon$ is a
Stein manifold, there exists an analytic right inverse operator
$R(z)$. Now, given $\phi (z) \in {\cal H}(\Phi_\epsilon)$, the
function $g(z)=R(z)\phi(z)$ satisfies $\tilde{t}_g=\phi$. This
proves the surjectivity that we need.

The last statement of the lemma about the mapping $t$ being a
topological homomorphism is just the open mapping theorem. \qed

\noindent \mbox{{\bf Proof of Lemma \ref{P-WC}}: } 1. We first
show that the operator $\cal{U}$ maps continuously the space $C_m$
into $C^{\infty}(T,H^m(K))$. Indeed, if $\varphi \in C_m$, then
$||\varphi||_{H^m(K+\gamma)}$ decays faster than any power of
$|\gamma|$. This together with (\ref{Floquet}) leads to the
immediate conclusion that $\cal{U} \varphi$ belongs to
$C^{\infty}(T,H^m(K))$ and to the continuity of the corresponding
mapping. Since $\cal{U} \varphi$ is a section of the sub-bundle
${\cal E}^m$ (see the Section 2.2 in \cite{Ku}), this gives us the
needed conclusion. Conversely, let $$ s(z) \in C^{\infty}(T,{\cal
E}^m) \subset C^{\infty}(T,H^m(K)) . $$ One can expand the
$H^m(K)$-valued function $s(z)$ into the Fourier series: $$
s(z)=\sum_{\gamma \in \Gamma }s_{\gamma}z^\gamma ,z\in T. $$ Here
$s_{\gamma} \in H^m(K)$. Standard estimates of the Fourier
coefficients of smooth functions apply, which show that
$||s_{\gamma}||$ decays faster than any power of $|\gamma|$. Let
us define now a function $\phi$ on $\Real^n$ such that
$\phi(x-\gamma)= s_{\gamma}(x)$ for $x \in K$ and $\gamma \in
\Gamma$. The additional information that $s$ is a section of the
sub-bundle ${\cal E}^m$ leads (as in \cite[page 96]{Ku}) to the
conclusion that $\phi \in H^m_{loc}(\Real^n)$. This implies that
$\phi \in C_m$ and finishes the proof of the first statement of
the lemma.

Statements 2 and 3 are correspondingly parts of Theorem 2.2.3 and
3.1.5 of \cite{Ku}. \qed

\noindent \mbox{{\bf Proof of Lemma \ref{delta}}: } The first
statement is rather obvious. Indeed, the statement is local, and
locally one can construct a smooth one-sided inverse. The second
statement can be proven like the similar statement in
\cite[Corollary 1.7.2]{Ku}. For completeness, we provide the
scheme of the proof here. Under the conditions of the second
statement of the lemma, it is easy to see that any functional
annihilating the range of the operator of multiplication by ${\cal
P}(z)$ must be supported at the finite set ${\cal Z}$ where ${\cal
P}(z)$ is not surjective. This also reduces the considerations to
a neighborhood $U$ of a point $z_0\in {\cal Z}$. Using the
Fredholm property, one can find a closed subspace $M$ of finite
codimension in $B_1$ such that the operators ${\cal P}(z)$ have
zero kernel on $M$ for all $z \in U$ (see the corresponding lemma
in \cite{At}, or Lemma 1.2.11 and Remark 2 below it in \cite{Ku}).
Now the problem reduces to a similar one on a finite-dimensional
space, where a standard representation of distributions supported
at a point implies (\ref{repr}). \qed

\mysection{Further remarks} \label{remarks}
 \begin{Rems}\label{Rems}
{\em 1.
 Throughout the paper, we have assumed for simplicity that
all the coefficients of the operators $P$ and $P^{*}$ are
$C^\infty$-smooth. In fact, we do not need such a restrictive
assumption (see the discussion in \cite[Section 3.4.D]{Ku}). For
example, a sufficient (but not necessary) condition for all the
statements of Section \ref{Sec3} to hold true is that the
coefficients of $L$ and $L^{*}$ are H\"{o}lder continuous.
Actually, even less is needed. For instance, conditions imposed on
the Schr\"{o}dinger operators in Theorem \ref{applicat} are
sufficient. It is clear that the conditions on the coefficients
could be significantly relaxed, if the operators were considered in
the weak sense, or by means of their quadratic forms. This should
not change the general techniques of the proofs.
We did not intend, however, to find the optimal
requirements on the coefficients for all our results to hold.
\vspace{2mm}

\noindent 2. It should be possible to describe the class of solutions
of the equation $Lu=0$ that are representable by a distribution
rather than by a hyperfunction. We plan to address this
problem elsewhere. \vspace{2mm}

\noindent 3. The Liouville theorem can probably be extended to
systems of equations (for instance, to the Maxwell system). In
this case one would face the problems of a possibly
nonzero index of the
corresponding operator and of multiple eigenvalues (the latter
can also occur for scalar operators). We believe
that the technique of this paper might be adjusted to handle
some of these situations. The extensions of the result of
\cite{KS} to the Pauli
and Maxwell operators obtained in \cite{BS1} and \cite{BS2} would
provide examples where the needed information on the behavior of
the dispersion relations at the bottom of the spectrum is
available. }
\end{Rems}
\appendix
 \mysection{Appendix} \label{Appendix2}
 In this appendix we present an alternative
proof of the third statement of Theorem \ref{applicat} in the case
when either $\Lambda_0=0$ and $N\geq 0$, or $\Lambda_0>0$ and
$0\leq N\leq 1$. The proof relies on some basic notions of
homogenization theory \cite{JKO} and imitates the proof of
Theorem 2 in \cite{MS}, where $L$ is assumed to be an operator in
divergence form. Therefore, we skip some details which are
essentially the same as in \cite{MS}.

We need to recall some basic definitions from homogenization
theory (see, for example, \cite{BLP,JKO}). Suppose that $L$ is a
second order elliptic operator of the form
\begin{equation}
L=-\sum_{i,j=1}^n a_{ij}(x)\partial _{i}\partial _{j}+\sum_{i=1}^n
b_i(x)\partial _i\, , \label{Aoperator}
\end{equation}
with periodic coefficients and denote the positive matrix
$\{a_{ij}(x)\}$ by $\mathcal{A}(x)$ and the periodic vector
$(b_1,\ldots,b_n)^T$ by $b$. Let $\psi$ be the positive normalized
periodic solution of the equation $L^*u=0$. Let
$\Psi(x)=(\Psi_1(x)\ldots,\Psi_n(x))^T$ be a solution of the
equation
\begin{equation}\label{Psi1}
L\Psi=-b(x)+\int_{\mathbb{T}^n}b(x)\psi(x)\,dx\quad \mbox{ in }
\mathbb{T}^n.
\end{equation}
Consider the matrix
\begin{equation}\label{Qeq}
  \mathcal{Q}=\{q_{ij}\}:=\int_{\mathbb{T}^n}(I+\nabla
\Psi)^T\mathcal{A}(x)(I+\nabla\Psi)\psi(x)\,dx\,,
\end{equation}
were $I$ is the identity matrix. The operator $Q:=-\sum_{i,j=1}^n
q_{ij}\partial _{i}\partial _{j}$ is called the {\em homogenized
operator} of the operator $L$, and the positive matrix
$\mathcal{Q}=\{q_{ij}\}$ is called the {\em homogenized matrix}
(see, \cite[Section 2.5]{JKO}).

The following lemma, which is actually a new formulation of
\cite[Theorem 5]{Pr}), establishes a connection between the
function $\Gl$ and homogenization theory.
\begin{Lem}\label{lemhomogenization}
Let $L$ be an operator of the form (\ref{operator}) and suppose
that $\xi\in\Xi$. Let $u_\xi$ and $u^*_{-\xi}$ be the positive
Bloch solutions of the equations $Lu=0$ and $L^*u=0$,
respectively. Denote by $\psi$ the periodic function $u_\xi
u^*_{-\xi}$. Consider the operator
\begin{equation}\label{operatortilde}
\tilde{L}=(u_\xi(x))^{-1}Lu_\xi(x)=-\sum_{i,j=1}^n
a_{ij}(x)\partial _{i}\partial _{j}+\sum_{i=1}^n
\tilde{b}_i(x)\partial _i\, ,
\end{equation}
let
\begin{equation}
Q=-\sum_{i,j=1}^n q_{ij}\partial _{i}\partial
_{j}\label{operatorhomo}
\end{equation}
be the homogenized operator of the operator $\tilde{L}$, and
$\mathcal{Q}=\{q_{ij}\}$ be the homogenized matrix. Then $\psi$ is
the principal eigenfunction of the operator $\tilde{L}^*$ on the
torus $\mathbb{T}^n$ with an eigenvalue $0$. Moreover,
$\Hess(\Lambda(\xi))=-\mathcal{Q}$.
\end{Lem}
\pf The first statement of the lemma can be checked easily while
the second statement follows directly from the formula in
\cite[Theorem 5]{Pr}, and the definition of the homogenized
operator. \qed

\vskip 3mm \noindent \mbox{{\bf Proof of a part of
the third statement of Theorem \ref{applicat}}: }
We clearly may assume that $L{\bf 1}=0$, so,
$$L=-\sum_{i,j=1}^n a_{ij}(x)\partial _{i}\partial
_{j}+\sum_{i=1}^n b_i(x)\partial _i .$$ We denote by  $\psi$ the
normalized positive solution of the equation $L^*u=0$ in
$\mathbb{T}^n$. Let $\Psi$ be a solution of the system
(\ref{Psi1}), and $Q$ be the homogenized operator of the operator
$L$.

Assume first that $\Gl_0\geq 0$. The case $N=0$ is trivial, and
follows from Theorem \ref{Liouville} and Lemma \ref{constants}.
Let $N=1$. Recall that according to Theorem \ref{Liouville},
$d_1\leq n+1$.
Moreover, by Theorem \ref{Liouville} and the Leibnitz's rule, a
(real) solution of linear growth is of the form
$$u(x)=\sum_{j=1}^n a_jx_j+\gf(x),$$ where $a_j\in \Real$ and
$\gf$ is periodic.

By lemma \ref{Lambdazero}, $\Gl_0=0$  if and only if for every
$1\leq k\leq n$
\begin{equation} \label{Pinsky}
\ga_j:=\int\limits_{\mathbb{T}^n}b_j(x)\psi(x)\,dx=0\,.
\end{equation}
For $1\leq j\leq n$, we write an ``Ansatz'' for a solution of
linear growth of the form
\begin{equation}\label{Fkeq}
F_j(x)=x_j+\gf_j(x),
\end{equation}
 where $\gf_j$ is a periodic function.
Clearly, $F_j$ is a solution of $Lu=0$ in $\mathbb{R}^n$ if and
only if $\gf_j(x)$ solves the nonhomogeneous equation $Lu=-b_j$ in
$\mathbb{T}^n$. By the Fredholm alternative, this equation is
solvable in $\mathbb{T}^n$ if and only if $\ga_j=0$ which holds
true for {\em all} $1\leq j\leq n$, if and only if $\Lambda_0=0$
(and in this case, $\gf_j=\Psi_j$, see (\ref{Psi1})). Therefore,
$d_1=n+1$ if $\Lambda_0=0$, and $d_1<n+1$ if $\Lambda_0>0$.

In order to finish the proof for $N=1$, we need to prove that if
$\Lambda_0>0$, then $d_1\geq n$. Without loss of generality, we
may assume that $\ga_n\neq 0$. We construct $(n-1)$ linearly
independent solutions of linear growth of the form
 $$F_j(x)=x_j-\ga_j(\ga_n)^{-1}x_n+\gf_j(x),$$
where $1\leq j\leq n-1$, and $\gf_j$ solves the equation
$Lu=-b_j+\ga_j(\ga_n)^{-1}b_n$. Note that these $(n-1)$ equations
are solvable and therefore, $d_1\geq n$.

For $N\geq 2$, we assume that $\Lambda_0=0$. Recall that if $u\in
V_N$ then by Theorem \ref{Liouville} and the Leibnitz's rule
$$u(x)=u^{(N)}(x) + \sum\limits_{|\gn |< N}x^\gn p_\gn (x) ,$$
where $$u^{(N)}(x)=\sum\limits_{|\gn |= N}x^\gn p_\gn \, ,$$ and
$p_\gn$ are periodic functions if $|\gn|<N$, and $p_\gn
\in\mathbb{R}$, if $|\gn |= N$.

\vskip 2mm \noindent \mbox{{\bf Claim}:} {\em Assume that
$\Lambda_0=0$. Then for all $N\geq0$
\begin{equation}\label{claimeq}
Qu^{(N)}=0\,.
\end{equation}
In particular,  $d_N\leq h_{n,N}$.}

\vskip 2mm

\noindent \mbox{{\bf Proof of the claim}:} Assume first that
$N=2$. Then $u\in V_2$ is of the form $$u(x)=\frac{1}{2}(Cx\cdot
x)+\sum_{j=1}^n x_jp_j(x)+p_0(x),$$ where $C$ is a constant
symmetric matrix, and $p_0,p_1,\ldots,p_n$ are periodic functions.

A direct calculation shows that the vector $p=(p_1,\ldots,p_n)^T$
must satisfy the equation $Lp=-Cb$ which is solvable since
$\Gl_0=0$. Therefore, $p=C\Psi$ (up to a constant vector). Also,
$p_0$ must satisfy
 $$Lp_0=f:=tr(A(I+2\nabla \Psi^T)C^T)-b\cdot
C\Psi\,.$$
 The compatibility condition for this equation is
$\int\limits_{\mathbb{T}^n}f(x)\psi(x)\,dx=0$ which after some
calculations implies that $$tr(\mathcal{Q}C^T)=0\,,$$ where
$\mathcal{Q}$ is the homogenized matrix of the operator $L$ (see
(\ref{Qeq})). Since $u^{(2)}:=\frac{1}{2}(Cx\cdot x)$ is a
homogeneous polynomial of degree $2$, it follows that $Qu^{(2)}=
tr(\mathcal{Q}C^T)$. Therefore, $u^{(2)}$ solves the equation
$Qu=0$. Thus, the case $N=2$ is settled.

For $N>2$, we proceed by induction as in \cite{MS}. Namely, assume
that the claim (\ref{claimeq}) has been proven for $N-1$, and let
$u\in V_N$. Let $\Delta_i$ be the difference operator $\Delta_i
f(x):=f(x+e_i)-f(x)$, where $e_i$ is the $i$-th vector of the
standard basis of $\mathbb{R}^n$, and $1\leq i\leq n$. Then
$v_i:=\Delta_i u\in V_{N-1}$ and the leading part of $v_i$ is
given by $(\Delta_i u)^{(N-1)}=\partial_i u^{(N)}$. By the
induction hypothesis, $Q((\Delta_i u)^{(N-1)})=0$. Therefore,
$$\partial_i(Qu^{(N)})=Q(\partial_i u^{(N)})=Q((\Delta_i
u)^{(N-1)})=0\quad 1\leq i\leq n\,.$$ Hence,
$Qu^{(N)}=\mbox{const.}\;$, and since $Qu^{(N)}$ is homogeneous of
degree $N-2>0$, we obtain that $Qu^{(N)}=0$, and the claim is
proved.

It remains to prove that $d_N\geq h_{n,N}$. So, for any
homogeneous polynomial $h$ of degree $N$ which is $Q$-harmonic, we
need to find a solution $u\in V_N$ such that $u^{(N)}=h$. Let
$u\in V_N$ and $\varepsilon>0$. Consider the function
 $$\vge^Nu(\frac{x}{\vge})=\sum_{|\gn|\leq N}
 \vge^{N-|\gn|}x^\gn p_\gn(\frac{x}{\vge}),$$
which tends to $u^{(N)}$ as $\vge\to 0$. We consider $x$ and
$y=\frac{x}{\vge}$ as independent variables and write
$$U(x,y,\vge):=\sum_{|\gn|\leq N}
 \vge^{N-|\gn|}x^\gn p_\gn(y)=U_0(x)+\vge U_1(x,y)+\cdots+
 \vge^NU_N(x,y)\,.$$
Then the equation $L(x,\partial_x)u=0$ implies that
 $$(L_0+\vge L_1+\vge^2 L_2)U=0\,,$$ where
 $$ L_0\!=\!L(y,\partial_y)\;;\,
L_1\!=\!-2\!\sum_{i,j=1}^n\!a_{ij}(y)\partial^2 _{x_{i},y_{j}}+
 \!\sum_{i=1}^nb_i(y)\partial_{x_{i}}\;;\,
L_2\!=\!-\!\!\sum_{i,j=1}^n\! a_{ij}(y)\partial^2_{x_{i},x_{j}}.$$

We look for a formal differential operator
 $$\Phi=\sum_{j=0}^\infty \vge^{k}\Phi_j=
 \sum_\gn \vge^{|\gn|}\phi_\gn(y)\partial_x^\gn\,,$$
 where $\phi_\gn(y)$ are periodic functions and $\phi_0=1$. This operator should
 satisfy
\begin{equation}\label{Lpsi}
(L_0+\vge L_1+\vge^2 L_2)\Phi=M+L_0(y,\partial_y) - 2\vge
\sum_{i,j=1}^n a_{ij}(y)\partial _{x_{i}}\partial _{y_{j}}\,,
\end{equation}
where the formal operator $$M=\sum_{j=2}^\infty
\vge^{j}M_j=\sum_{|\gn|\geq 2}
\vge^{|\gn|}m_\gn\partial_x^\gn\,,$$ has constant coefficients.

Comparing the coefficients of $\vge^s$ in (\ref{Lpsi}) yields the
following equations (the equation for $s=0$ is automatically
satisfied).
 \bea\label{1eq}
 L_0\Phi_1+L_1= -2
\sum_{i,j=1}^n a_{ij}(y)\partial _{x_{i}}\partial _{y_{j}}
\,,\quad & s=1\,,&\\
 L_0\Phi_s+L_1\Phi_{s-1}+L_2\Phi_{s-2}=M_s\,,\quad & s\geq 2\,.&
 \label{seqs}\eea
It is easily checked that for $s=1$ the functions $\phi_j(y)$ of
Equation (\ref{Fkeq}) are the corresponding solutions for
$\Phi_1$. Also, Equation (\ref{seqs}) for $s=2$ is solvable if
$M_2=Q$, where $Q$ is the homogenized operator of $L$. Similarly,
the constant coefficients of the operator $M_s\,,\; s>2$, are
determined by the compatibility condition for Equation
(\ref{seqs}) with $s>2$.

Let $R\,:\, \mathcal{P}\to \mathcal{P}$ be a linear right inverse
of the homogenized operator $Q$ that preserves the homogeneity of
polynomials. Consider the formal operator $A$ which is defined by
the equation
 $$A-I=R\sum_{j=1}^\infty \vge^j M_{j+2} \,,$$
 and let $A^{-1}$ be its unique formal inverse. Note that $\vge^2 M_2A=M$.

Let $U_0(x)$ be a given homogeneous polynomial of degree $N$ which
solves the equation $Qu=0$, and let $V(x):=A^{-1}U_0(x)$. We have
 $$MV=\vge^2M_2AV=\vge^2M_2U_0=0.$$
Define $U(x,y,\vge):=\Phi A^{-1}U_0=\Phi V$, and denote
$u(x):=U(x,x,1)$. By inspection, $u$ has a polynomial growth of
order $N$, and $u^{(N)}(x)=U_0(x)$. Moreover, \beqanl (L_0+\vge
L_1+\vge^2 L_2)U=(L_0+\vge L_1+\vge^2 L_2)\Phi
V\\=MV+L_0(y,\partial_y)V(x) + 2\vge \sum_{i,j=1}^n
a_{ij}(y)\partial _{x_{i}}\partial _{y_{j}}V(x)=0\,, \eeqanl
 and $\Phi A^{-1}$ is the desired mapping.
 \qed

\begin{Rem}\label{Ancona}
{\em 1. Let $F_j$ be the solutions of linear growth defined by
Equation (\ref{Fkeq}). A.~Ancona \cite{An} proved that the map
$F(x)=(F_1(x),\dots,F_n(x))$ is a diffeomorphism on $\Real^n$ if
$n\leq 2$, while for $n>2$ this map is not necessarily a
diffeomorphism. \vspace{2mm}

\noindent 2. Assume that $L{\bf 1}=0$ and $\Lambda_0=0$. Let
$\Lambda(\xi)= \sum_{|\gn|\geq 2}a_\gn\xi^\gn$ be the Taylor
expansion of the function $\Lambda$. We conjecture that
$a_\gn=m_\gn$, where $m_\gn$ are the coefficients of the operator
$M$.
 }\end{Rem}
\begin{center}
{\bf Acknowledgments} \end{center} The authors express their
gratitude to Professors S. Agmon and V. Lin for useful discussions
and to Professor P.~Li for the information about the manuscript
\cite{LW1}.

The work of P. Kuchment was partially supported by the NSF Grant
DMS 9610444 and by a DEPSCoR Grant. P. Kuchment expresses his
gratitude to NSF, ARO, and to the State of Kansas for this
support. The content of this paper does not necessarily reflect
the position or the policy of the federal government of the USA,
and no official endorsement should be inferred. The work of
Y.~Pinchover was partially supported by the Fund for the Promotion
of Research at the Technion.

\end{document}